\documentclass[a4paper,10pt]{elsarticle}
\usepackage{framed,multirow}
\usepackage[utf8]{inputenc}
\usepackage{bm}
\usepackage{amsmath}
\usepackage{mathrsfs}
\usepackage{graphicx}
\usepackage{epsfig, setspace}
\usepackage{url}
\usepackage{graphicx, transparent, color}
\usepackage{algorithm} 
\usepackage{algorithmicx}
\usepackage{etex}
\usepackage[bookmarks=false]{hyperref}
\usepackage{amsmath}
\usepackage{xcolor}
\usepackage{rotating} 
\usepackage{pdflscape} 
\usepackage{silence}
\usepackage{ulem,cancel}
\usepackage[utf8]{inputenc}
\newtheorem{remark}{\bf Remark}[section]
\newtheorem{example}{Example}[section]
\usepackage{graphicx}
\usepackage{float}
\usepackage{subfigure}
\usepackage{subfigmat}

\usepackage{tikz}
\usepackage{capt-of}
\usepackage{setspace,lipsum}
\newcommand{\half}{\frac{1}{2}}

\usepackage[top=18truemm,bottom=18truemm,left=20truemm,right=20truemm]{geometry}

\begin{document}
\begin{frontmatter}
\title{\textcolor{black}{Implicit gradients based conservative numerical scheme for compressible flows}}
\author[AA_address]{Amareshwara Sainadh Chamarthi \cortext[cor1]{Corresponding author. \\ 
E-mail address: skywayman9@gmail.com (Amareshwara Sainadh  Ch.).}}
\author[AA_address]{Natan\ Hoffmann}
\author[CC_address]{Hiroaki Nishikawa}
\author[AA_address]{Steven H.\ Frankel}
\address[AA_address]{Faculty of Mechanical Engineering, Technion - Israel Institute of Technology, Haifa, Israel}
\address[CC_address]{National Institute of Aerospace, 100 Exploration Way, Hampton, VA 23666, USA}

\begin{abstract}

\textcolor{black}{This paper introduces a novel approach to compute the numerical fluxes at the cell boundaries for a cell-centered conservative numerical scheme. Explicit gradients used in deriving the reconstruction polynomials are replaced by high-order gradients computed by compact finite differences, referred to as implicit gradients in this paper. A problem-independent shock capturing approach via Boundary Variation Diminishing (BVD) algorithm is used to suppress oscillations for the simulation of flows with shocks and material interfaces. Several numerical test cases are carried out to verify the proposed  method's capability using the implicit gradient method for compressible flows.}

\end{abstract}


\end{frontmatter}
\section{Introduction}\label{sec-1}
 
This paper presents a novel algorithm based on implicit gradients for regular hexahedral (i.e., Cartesian) grids to ensure low dispersion and dissipation. {\color{black} It is demonstrated that such low dispersion and dissipation schemes can be constructed by {\it linearly high-order} schemes that achieve high-order accuracy for linear equations but second-order accurate at best for nonlinear equations. Specifically, we construct low dispersion and dissipation schemes based on a finite-volume scheme with the kappa solution reconstruction scheme of Van Leer \cite{van1977towards,van1985upwind} by making two key modifications: (1) replace cell-averaged solutions by point-valued solutions and (2) express the kappa scheme with first- and second-derivatives and replace them by those computed by implicit methods. The first modification is introduced to construct a high-order conservative finite-difference scheme, rather than a finite-volume scheme, so that high-order accuracy can be achieved (for linear equations) in multi-dimensions simply by applying a one-dimensional scheme in each coordinate direction. Implicit gradient schemes are similar to that of compact finite volume schemes proposed by Sengupta et al. \cite{sengupta2005new}. It is similar in that it uses the kappa-reconstruction scheme, Equation (\ref{eqn:3linear}) in this paper, with implicit gradients (see Sengupta et al. \cite{sengupta2005new} Equation 1). The difference is that  Sengupta et al. \cite{sengupta2005new} have used upwind compact finite difference schemes to compute the first derivatives, whereas, in this paper, the first derivatives are computed by implicit central finite difference schemes. As pointed out in Refs \cite{VanLeerNishikawa_UltimateUnderstanding:JCP2021,FalseAccuracyUMUSCL:CiCP2021}, however, methods based on the finite-volume method with point-valued solutions (instead of cell-averaged solutions) can achieve third- or higher-order accuracy for linear equations but can only be second-order accurate for nonlinear equations. While high-order accuracy can be achieved for nonlinear equations by introducing high-order flux reconstruction \cite{VanLeerNishikawa_UltimateUnderstanding:JCP2021,FalseAccuracyUMUSCL:CiCP2021,Nishikawa_FSR:2020}, the linearly high-order schemes are obviously more computationally efficient since expensive flux reconstruction is not necessary. Moreover, it has been demonstrated that linearly high-order schemes significantly improve the resolution of complex three-dimensional turbulent-flow simulations at a little additional cost over conventional second-order finite-volume schemes widely used in practical unstructured-grid computational fluid dynamics codes, despite not being genuinely high-order accurate for such nonlinear problems \cite{yang_harris:AIAAJ2016,HQYANG:AIAA2013-2021,yang_harris:CCP2018,GarciaBarakos:IJNMF2018,zhang2011order,DementRuffin:aiaa2018-1305}. }

The proposed method is conservative and based on an upwind flux computed with solutions reconstructed by the kappa scheme of Van Leer \cite{van1977towards,van1985upwind}. {\color{black} This method is expected to be third-order accurate because the kappa scheme is a quadratic reconstruction scheme. However, contrary to expectations, it results in a {\it fourth-order} upwind finite-difference scheme {\color{black} (at least for linear equations)} if combined with high-order accurate gradients on regular grids. Specifically, we will express the kappa scheme in terms of the first and second derivatives of the solution and then compute these derivatives by high-order gradient methods. In particular, conservative schemes with exceptionally low dispersion and dissipative errors can be obtained if the gradients are computed implicitly (globally coupled linear systems). Implicit gradient methods are not new: high-order compact schemes proposed by \textcolor{black}{Lele} \cite{lele1992compact} can be directly applied to compute the gradients on regular grids. {\color{black} However, in our schemes, implicit gradients are used not for directly approximating the flux divergence but for evaluating the derivatives in solution reconstruction schemes within a framework of conservative \textcolor{black}{finite-difference-type} schemes. It allows us to easily construct stable high-resolution shock-capturing schemes incorporating Riemann solvers and various monotonicity-preserving mechanisms in the solution reconstruction, which would not be simple to incorporate if the flux divergence was directly computed by implicit gradient methods as typical in the so-called compact schemes.} In this paper, we will employ the fourth- and sixth-order compact finite difference schemes of Lele \cite{lele1992compact} as implicit gradient methods. More specifically, we define the numerical solutions as point values at cell centers (not cell averages) and evaluate the numerical flux with solution values reconstructed at a face by a quadratic Legendre polynomial {\color{black} as in} the unlimited kappa-scheme of Van Leer \cite{van1977towards,van1985upwind} with first and second derivatives computed by implicit gradient methods \cite{lele1992compact}. For shock capturing, we will combine the proposed schemes with the BVD algorithm \cite{sun2016boundary}. The proposed method has the following advantages:

\begin{description}
	\item[(a)] It generates fourth-order upwind {\color{black} finite-difference} schemes with a quadratic reconstruction {\color{black} for linear equations}, which typically leads to third-order accuracy at best. It generates practical low-dispersion/dissipation schemes that can be easily implemented for structured-grid codes, and
	\item[(b)] the implicit gradient approach combined with the shock-capturing approach via the Boundary Variation Diminishing (BVD) algorithm  {\color{black} gives} superior results of flows with shocks, material interfaces and small scale features than the approach presented in \cite{chamarthi2021high}.
\end{description}

{\color{black} The BVD algorithm was initially proposed by Sun et al. \cite{sun2016boundary} which combines a non-polynomial reconstruction scheme, THINC (Tangent of Hyperbola for INterface Capturing), for discontinuous regions and an unlimited polynomial based reconstruction for the smooth regions of flows. The proposed methodology adaptively chooses the scheme with minimum Total Boundary Variation (TBV), reducing the numerical dissipation. Following their idea, Chamarthi and Frankel \cite{chamarthi2021high} presented a new algorithm named HOCUS (High-Order Central Upwind Scheme), which combined the \textcolor {black}{Monotonicity preserving (MP)} scheme and a linear-compact scheme using the BVD principle, which is used in the present approach.

The rest of the paper is organized as follows. Section \ref{sec-2} introduces the governing equations of viscous compressible flows. In Section \ref{sec-3}, a brief description of the \textcolor{black}{cell-centered conservative approach} is presented and {\color{black} the novel reconstruction schemes are introduced along with the implementation details}. Numerical results and discussion are presented in Section \ref{sec-4}, and finally, in Section \ref{sec-6}, we provide concluding remarks.

\section{Governing equations}\label{sec-2}

The compressible Navier–Stokes (NS) equations in a Cartesian coordinate system can be expressed as:

\begin{equation}\label{eq:scalar_conservation}
\frac{\partial \mathbf{Q}}{\partial t}+\frac{\partial \mathbf{F^c}}{\partial x}+\frac{\partial \mathbf{G^c}}{\partial y}+\frac{\partial \mathbf{H^c}}{\partial z}+\frac{\partial \mathbf{F^v}}{\partial x}+\frac{\partial \mathbf{G^v}}{\partial y}+\frac{\partial \mathbf{H^v}}{\partial z}=0,
\end{equation}
where $\mathbf{Q} = (\rho, \rho u, \rho v, \rho w, \rho E)^{T}$ is the conserved variable vector. $\rho$, $u$, $v$, $w$, and $\rho E$ are the density, the three Cartesian velocity components and the total energy, respectively. $E = e + \frac{1}{2}(u^2 + v^2 + w^2)$ is the total energy per unit mass, where $e$ is the internal energy per unit mass. $\mathbf{F^c}$, $\mathbf{G^c}$,  $\mathbf{H^c}$ and $\mathbf{F^v}$, $\mathbf{G^v}$, $\mathbf{H^v}$, are the convective (superscript $c$) and viscous (superscript $v$) flux vectors in each coordinate direction, respectively. The convective and viscous flux vectors are given as:

\begin{equation}
\begin{array}{l}\label{eqn-inv}
\mathbf{F^c}=\left[\rho u, \rho u^{2}+p, \rho u v , \rho u w, u(\rho E+p)\right]^{T}, \\
\mathbf{G^c}=\left[\rho v, \rho u v,\rho v^{2}+p, \rho v w, v(\rho E+p)\right]^{T}, \\
\mathbf{H^c}=\left[\rho w, \rho u w, \rho v w,\rho w^{2}+p, w(\rho E+p)\right]^{T},
\end{array}
\end{equation}

\begin{equation}
\begin{array}{l}\label{eqn-visc}
\mathbf{F^v}=\left[0, \tau_{x x}, \tau_{x y}, \tau_{x z}, u \tau_{x x}+v \tau_{x y}+w \tau_{x z}-q_{x}\right]^{T}, \\
\mathbf{G^v}=\left[0, \tau_{x y}, \tau_{y y}, \tau_{y z}, u \tau_{y x}+v \tau_{y y}+w \tau_{y z}-q_{y}\right]^{T}, \\
\mathbf{H^v}=\left[0, \tau_{x z}, \tau_{y z}, \tau_{z z}, u \tau_{z x}+v \tau_{z y}+w \tau_{z z}-q_{z}\right]^{T},
\end{array}
\end{equation}
The system of equations is closed with the ideal gas equation of state which relates the thermodynamic pressure $p$ and the total energy per unit mass:}
\begin{equation}\label{eqn:pressure}
p = (\gamma -1) (E - \rho \frac{(u^2+v^2+w^2)}{2}),
\end{equation}
where $\gamma$ is the ratio of specific heats of the fluid ($\gamma=1.4$ for air at standard conditions). After non-dimensionalizing the governing equations by reference quantities (for velocity we use the freestream velocity as the reference), the components of the viscous stress tensor $\tau$ and the heat flux $q$ are defined in tensor notation as:

\begin{equation}\label{eqn:5-stress}
\tau_{i j}=\frac{\mu}{\operatorname{Re}}\left(\frac{\partial u_{i}}{\partial x_{j}}+\frac{\partial u_{j}}{\partial x_{i}}-\frac{2}{3} \frac{\partial u_{k}}{\partial x_{k}} \delta_{i j}\right),
\end{equation}

\begin{equation}\label{eqn:6-heat}
\begin{aligned}
\mathrm{q}_{i}=-\frac{\mu}{\operatorname{Re Pr Ma}(\gamma-1)} \frac{\partial T}{\partial x_{i}},
\end{aligned}
\end{equation}

\noindent where $\mu$ is the dynamic viscosity, $\delta_{ij}$ is the Kronecker delta, Ma is the Mach number, Re is the Reynolds number, Pr is the Prandtl number, and $T$ is the temperature obtained from the ideal gas assumption:

\begin{equation}
    T= \text{Ma}^{2} \gamma \frac{p}{\rho}.
\end{equation}

\section{Numerical methods}\label{sec-3}

The time evolution of the vector of cell-centered conservative variables $\mathbf{\hat Q}$ is given by the following semi-discrete relation {\color{black} applied to a Cartesian cell $I_{j,i,k} = [x_{j-1/2}, x_{j+1/2}] \times [y_{i-1/2}, y_{i+1/2}] \times [z_{k-1/2}, z_{k+1/2}]$}, expressed as an ordinary differential equation:

\textcolor{black}{
 \begin{equation}\label{eqn-differencing}
\frac{\mathrm{d}}{\mathrm{dt}} {\mathbf{\hat Q}}_{i,j,k} = \mathbf{Res}_{i,j,k} = - \frac{d \mathbf F}{dx}_{i,j,k} - \frac{d \mathbf G}{dy}_{i,j,k}- \frac{d \mathbf H}{dz}_{i,j,k},
 \end{equation}}
 
\textcolor{black}{ where $ \mathbf{Res}_{i,j,k}$ denotes the residual. $\frac{d \mathbf F}{dx}_{i,j,k}$, $\frac{d \mathbf G}{dy}_{i,j,k}$ and $\frac{d \mathbf H}{dz}_{i,j,k}$ are approximations to the flux derivatives at the cell center, and we seek their high-order approximations in the conservative form:}
 \textcolor{black}{
 \begin{equation}\label{eqn-differencing_residual}
\begin{aligned}
\frac{d \mathbf F}{dx}_{i,j,k}  = &\frac{1}{\Delta x}\left[\left(\mathbf {\hat{F}^c}_{i+ \frac{1}{2}, j, k}-\mathbf {\hat{F}^c}_{i- \frac{1}{2}, j, k}\right) - \left(\mathbf {\hat{F}^v}_{i+ \frac{1}{2}, j, k}-\mathbf {\hat{F}^v}_{i-\frac{1}{2}, j, k}\right)\right] \\
\frac{d \mathbf G}{dy}_{i,j,k}=&\frac{1}{\Delta y}\left[\left(\mathbf {\hat{G}^c}_{i, j+ \frac{1}{2}, k}-\mathbf {\hat{G}^c}_{i, j- \frac{1}{2}, k}\right)-\left(\mathbf {\hat{G}^v}_{i, j+ \frac{1}{2}, k}-\mathbf {\hat{G}^v}_{i, j- \frac{1}{2}, k}\right)\right]\\
\frac{d \mathbf H}{dz}_{i,j,k}=&\frac{1}{\Delta z}\left[\left(\mathbf {\hat{H}^c}_{i, j, k+ \frac{1}{2}}-\mathbf {\hat{H}^c}_{i, j, k- \frac{1}{2}}\right)-\left(\mathbf {\hat{H}^v}_{i, j, k+ \frac{1}{2}}-\mathbf {\hat{H}^v}_{i, j, k- \frac{1}{2}}\right)\right],
\end{aligned}
\end{equation}}

where $\Delta x=x_{j+1 / 2}-x_{j-1 / 2}$, $\Delta y=y_{i+1 / 2}-y_{i-1 / 2}$, and $\Delta z=z_{k+1 / 2}-z_{k-1 / 2}$. $\mathbf {\hat{F}^c}$, $\mathbf {\hat{G}^c}$, $\mathbf {\hat{H}^c}$  and $\mathbf {\hat{F}^v}$, $\mathbf {\hat{G}^v}$ and $\mathbf {\hat{H}^v}$ are interpreted as the numerical approximation of the convective and viscous fluxes in the $x-$, $y-$ and $z-$directions, respectively. The viscous fluxes are computed using the $\alpha$-damping approach presented by \textcolor{black}{Chamarthi et al.} \cite{chamarthi2022}, which is extended in \cite{chamextend} for implicit gradients. This viscous flux discretization approach prevents odd-even decoupling and improves the solution quality of the viscous simulations.

In the following subsections, we provide the details of the computation of convective fluxes, including the novel implicit gradient method (\ref{sec-3.1}), shock-capturing algorithm (\ref{sec-3.2}), and the details of the approximate Riemann solver (\ref{sec-3.3}). {\color{black} Note that we define the cell-centered numerical solution $\mathbf{\hat Q}$ as a point-valued solution at the cell-center, not a cell-averaged solution. It simplifies the introduction of implicit gradients, the implementations of boundary conditions and initial values, {\color{black} and most importantly, it leads to linearly high-order finite-difference schemes that preserve the design order of accuracy in multi-dimensions. Hence, the fluxes are computed only at the center of a face in the $x$-direction as in a one-dimensional scheme.}

\subsection{Spatial discretization of fluxes}\label{sec-3.1}
 In this section, we present the spatial discretization of the inviscid fluxes. {\color{black} As mentioned, we desribe the method in one-dimension as it can be easily extended to multi-dimensional (2D and 3D) problems as a finite-difference scheme via dimension by dimension approach.} The governing equations (\ref{eq:scalar_conservation}) are discretized on a uniform grid with $N$ cells on a spatial domain spanning $x \in \left[x_a, x_b \right]$. The cell center locations are at $x_j = x_a + (j - 1/2) \Delta x$, $\forall j \in \{1, \: 2, \: \dots, \: N\}$, where $\Delta x = (x_b - x_a)/N$. The cell interfaces, indexed by half integer values, are at $x_{j+\frac{1}{2}}$, $\forall j \in \{ 0, \: 1, \: 2, \: \dots, \: N \}$. Let $I_j = [x_{j - 1 / 2}, x_{j + 1 / 2}]$ be a control volume (a computational cell) of width $\Delta x  = x_{j + 1 / 2} - x_{j - 1 / 2}$.}

 \subsubsection{Upwind flux (Riemann solver)}\label{sec-3.1.1} 
 
 The convective numerical fluxes in the Equation (\ref{eqn-differencing_residual}), $\mathbf {\hat{F^c}}_{j- 1 / 2}$ and $\mathbf {\hat{F^c}}_{j+ 1 / 2}$, are computed by an approximate Riemann solver. There are several types of Riemann solvers in the literature  \cite{deledicque2007exact, ivings1998riemann,roe1981approximate, batten1997choice, einfeldt1988godunov, toro1994restoration, osher1982upwind}, and the canonical form of the Riemann flux can be written as:

\begin{equation}\label{Numerical Flux}
\mathbf {\hat{F^c}}_{j+ 1 / 2}={F}^{\rm Riemann}_{j+\frac{1}{2}}(\mathbf{Q}_{j+\frac{1}{2}}^{L},\mathbf{Q}_{j+\frac{1}{2}}^{R}), 
\end{equation}

\begin{equation}
\mathbf{F}^{\rm Riemann}_{j+\frac{1}{2}}
= \frac{1}{2}
{\color{black}
\left[
{\mathbf{F}}(\mathbf{Q}_{j+\frac{1}{2}}^{L}) 
+ 
{\mathbf{F}}(\mathbf{Q}_{j+\frac{1}{2}}^{R})
\right]
}
-
 \frac{1}{2} | {\mathbf{A}_{j+\frac{1}{2}}}|({\mathbf{Q}^R_{j+\frac{1}{2}}}-{\mathbf{Q}^L_{j+\frac{1}{2}}}),
\label{eqn:Riemann}
\end{equation}
where $L$ and $R$ are adjacent values {\color{black} of a reconstructed solution polynomial at} a cell interface and $|{\mathbf{A}_{j+\frac{1}{2}}}|$ denotes the characteristic signal velocity evaluated at the cell interface in a hyperbolic equation. In our implementation, we obtain the cell interface conservative variable vector from the cell interface primitive variable vector,  $\mathbf{U}_{j+\frac{1}{2}}$ = $(\rho, u,  v, p)^T$. Thus, we present our reconstruction approach in terms of the primivitive variable vector.

\begin{remark}\label{eqn:whichgradeints}
\normalfont  \textcolor{black}{Both conservative, $\mathbf{Q}$, and primitive variables, $\mathbf{U}$, can be used for the evaluation of the gradients in the implicit gradient approach. For shock-capturing purposes, which will be explained in Section \ref{sec-3.1.2},  the gradients of the primitive variables are used for reconstruction.} 
\end{remark}

 \subsubsection{Reconstruction with the kappa scheme}\label{sec-3.1.2} 
 
The procedure of obtaining the values at the interface from \textcolor{black}{cell center variables} is called \textcolor{black}{reconstruction or interpolation}. It is obvious from Equation (\ref{eqn:Riemann}), a core problem is how to reconstruct the left- and right-side values, $\mathbf{Q}_{j+\frac{1}{2}}^{L}$ and $\mathbf{Q}_{j+\frac{1}{2}}^{R}$, for cell boundaries, which can fundamentally influence the numerical solution. 
Representing these numerical approximations of $L$ and $R$ at cell interface as a piecewise constant is equivalent to first-order 
approximation, i.e.,
\begin{equation}\label{eqn:first-order}
\begin{aligned}
{\mathbf{U}^L}_{j+\frac{1}{2}} &= \mathbf{\hat{U}_j}, \\
{\mathbf{U}^R}_{j+\frac{1}{2}} &= \mathbf{\hat{U}_{j+1}}. 
\end{aligned}
\end{equation}

Therefore, a linear approximation is a second-order spatial approximation, while a quadratic representation on each cell leads to a third-order spatial approximation. By considering a general local representation, as explained in \cite{van1977towards}, {\color{black} the quadratic} approximation can be expressed in terms of Legendre polynomials, valid for $x_{j-1 / 2} \leq x \leq x_{j+1 / 2}$:
\begin{equation}
 {\color{black} 
\mathbf{U}(x)=\mathbf{\hat{U}}_{j}+ {\mathbf{U}'_j} \left(x-x_{j}\right)+\frac{3  \mathbf{U}''_j}{2 } \kappa\left[\left(x-x_{j}\right)^{2}-\frac{\Delta x^{2}}{12}\right],
}
\label{quadratic_polynomial_at_j}
\end{equation}
where $\mathbf{\hat{U}}_{j}$ is the \textit{\textcolor{black}{cell-center value}} and $\mathbf{U}'_j$, ${\mathbf{U}''_j}$ are the estimations of the first and second derivatives within the cell $j$. Equation (\ref{quadratic_polynomial_at_j}) is the basis for the Monotonic Upstream-centered Scheme for Conservation Laws (MUSCL) scheme, popularly known as the kappa scheme of Van Leer \cite{van1977towards, van1985upwind}. The WENO schemes are an extension of the MUSCL scheme to an arbitrary order of accuracy (see derivations \textcolor{black}{\textcolor{black}{presented by Balsara et al. in}} \cite{balsara2016efficient}). {\color{black} Note, however, that the kappa scheme is used here with the point-values stored at the cell center, $\mathbf{\hat{U}}_{j}$, whereas it is used with the cell-average in MUSCL. Such a scheme can still achieve third-order accuracy for linear equations with $\kappa=\frac{1}{3}$, but can only be second-order accurate for nonlinear equations  \cite{VanLeerNishikawa_UltimateUnderstanding:JCP2021,FalseAccuracyUMUSCL:CiCP2021}. However, while it remains second-order accurate for nonlinear equations, it can achieve fourth-order accuracy for linear equations if the derivatives $\mathbf{U}'_j$ and $\mathbf{U}''_j$ are evaluated with high-order gradients. Linear fourth-order accuracy with a quadratic polynomial may be an unexpected result, however we provide a proof for this through Fourier analysis and multiple test cases.}

For the numerical approximations of the Riemann problem, we need the values at the cell interfaces only. Setting $x = x_j \pm \Delta x/2$ gives us the interface values:

\begin{equation}\label{eqn:left_right}
\begin{aligned}
\mathbf{U}_{j+1 / 2}^{L} 
&=\mathbf{\hat {U}}_{j}+\frac{\Delta x}{2} \mathbf{U}'_j+\frac{\kappa \Delta x^2}{4} \mathbf{U}''_j,  
{\color{black}
\quad \text{or} \quad  \ 
\mathbf{U}_{j-1 / 2}^{L} }
&=
{\color{black} \mathbf{\hat {U}}_{j-1} + \frac{\Delta x}{2} \mathbf{U}'_{j-1}+\frac{\kappa \Delta x^2}{4} \mathbf{U}''_{j-1} },
\\
\mathbf{U}_{j-1 / 2}^{R} &=\mathbf{\hat {U}}_{j}-\frac{\Delta x}{2}  \mathbf{U}'_j+\frac{\kappa \Delta x^2}{4}  \mathbf{U}''_j, \quad \text{or} \quad \ \mathbf{U}_{j+1 / 2}^{R} &=\mathbf{\hat {U}}_{j+1}-\frac{\Delta x}{2} \mathbf{U}'_{j+1}+\frac{\kappa \Delta x^2}{4} \mathbf{U}''_{j+1}.
\end{aligned}
\end{equation}
In order to define these approximations completely, the derivatives ${ \mathbf{U}'_j}$ and ${\mathbf{U}''_j}$ have to be estimated. {\color{black} Typically, these derivatives are computed via explicit finite differences. However, in this work, we compute the necessary derivatives using the implicit gradient approach. Both the explicit and implicit approaches will be presented in the next section.}

 \subsubsection{Second-order explicit gradients: Baseline scheme}\label{sec-3.1.3} 
 
 {\color{black} Before we present the implicit gradient approach, we describe for comparison the baseline MUSCL-type scheme, which uses second-order explicit finite differences for the gradients presented above.}  \textcolor{black}{By using $\kappa=\frac{1}{3}$ and substituting the following explicit central differences for the derivatives in Equation (\ref{eqn:left_right})},

\begin{equation}\label{eqn:expligrad}
\begin{aligned}
 \mathbf{U}'_j &= \frac{\mathbf{\hat {U}}_{j+1} - \mathbf{\hat {U}}_{j-1}}{2 \Delta x },\\
 \mathbf{U}''_j &=\frac{  \mathbf{\hat {U}}_{j+1} - 2\mathbf{\hat {U}}_{j} + \mathbf{\hat {U}}_{j-1} }{\Delta x^2} ,
\end{aligned}
\end{equation}
we obtain the following third order reconstruction {\color{black} formulas}, 

\begin{equation}
\begin{aligned} 
{\mathbf{U}^L}_{j+\frac{1}{2}}&= \frac{1}{6}\left(-\mathbf{\hat {U}}_{j-1} + 5\mathbf{\hat {U}}_{j} + 2\mathbf{\hat {U}}_{j+1} \right), \\
{\mathbf{U}^R}_{j-\frac{1}{2}} &= \frac{1}{6}\left(2\mathbf{\hat {U}}_{j-1} + 5\mathbf{\hat {U}}_{j} - \mathbf{\hat {U}}_{j+1} \right).
\end{aligned}
\label{eqn:3linear}
\end{equation}
{\color{black} These reconstructions lead to third-order accuracy for linear equations and second-order accuracy for non-linear equations because of the use of point values.}  Note that this should not be considered as a failure of the scheme design because it is just a simplified version of the original finite-difference scheme of Van Leer \cite{VanLeerNishikawa_UltimateUnderstanding:JCP2021} without flux reconstruction and it is deliberately designed this way in this work for developing practical and efficient schemes with low dispersion and dissipation. As we will show later, the scheme can achieve higher than third-order accuracy (for linear equations) with these explicit gradients replaced by higher-order implicit gradient formulas. {\color{black} In effect, the resulting schemes are linearly high-order conservative finite-difference schemes approximating the flux divergence in a conservative manner by a flux difference as shown in Equations (\ref{eqn-differencing_residual}), rather than the finite-volume scheme approximating the integral form of the governing equations}. Therefore, the quadratic reconstruction scheme does not limit the order of accuracy of the resulting scheme: it is used as a stepping stone towards a higher-order difference approximation. High-order explicit gradient formulas may be employed instead of implicit gradient methods, but such will significantly extend the residual stencil and require complicated algorithms near boundaries. Implicit gradients do not involve such complications near boundaries, and also have been found to generate much lower dispersion and dissipation that high-order explicit gradients.


 \subsubsection{High-order implicit gradients: Novel IG4 and IG6 schemes}\label{sec-3.1.5} 
 
{\color{black} Using the implicit gradient approach, we employ the compact schemes of Lele \cite{lele1992compact} to compute the gradients at cell centers.} For the first derivative, this can be written in general form as \cite{lele1992compact,nagarajan2003robust,boersma2005staggeblack}:

\begin{equation} \label{eqn:pade}\small
        \beta \mathbf{ {U}}'_{j-2} + \alpha \mathbf{ {U}}'_{j-1} + \mathbf{ {U}}'_{j} + \alpha \mathbf{ {U}}'_{j+1} + \beta \mathbf{ {U}}'_{j+2} =
        c \frac{\mathbf{\hat {U}}_{j+3} - \mathbf{\hat {U}}_{j-3}}{6 \Delta x} + b \frac{\mathbf{\hat {U}}_{j+2} - \mathbf{\hat {U}}_{j-2}}{4 \Delta x} + a \frac{\mathbf{ \hat{U}}_{j+1} - \mathbf{\hat {U}}_{j-1}}{2 \Delta x}.
\end{equation}
The left hand side of Equation (\ref{eqn:pade}) contains the spatial derivatives $\mathbf{{U}}'_j$ while the right hand side contains the function values $\mathbf{\hat {U}}$ at the cell center $x_j$. Compact finite difference schemes of different orders of accuracy are derived by matching the Taylor series coefficients with different constraints on the parameters $\alpha$, $\beta$, $a$, $b$ and $c$ and are listed in \cite{lele1992compact}. In this work, we considered difference schemes for the first derivatives with the following parameters,

\begin{subequations}
     \begin{alignat}{4}
&\beta = 0, \quad &a_1=\frac{2}{3}(\alpha+2), \quad &b_1=\frac{1}{3}(4 \alpha-1), \quad &c=0, \label{eq:param4}
     \end{alignat}
   \end{subequations}
By substituting $\alpha$ = $\frac{5}{14}$ in Equation (\ref{eq:param4}) we obtain the optimised fourth-order compact derivative (\textcolor{black}{see Fig. 2} in \cite{lele1992compact}), which we denote as \textbf{CD4} and which is written as: 

\begin{equation}
     \frac{5}{14} \mathbf{ {U}}_{j-1}^{\prime}+\mathbf{ {U}}_{j}^{\prime}+\frac{5}{14} \mathbf{ {U}}_{j+1}^{\prime}=\frac{b_1}{4 \Delta x}\left(\mathbf{\hat {U}}_{j+2}-\mathbf{\hat {U}}_{j-2}\right)+\frac{a_1}{2 \Delta x}\left(\mathbf{\hat {U}}_{j+1}-\mathbf{\hat {U}}_{j-1}\right). \label{eq:cd4}
\end{equation}

\noindent For $\alpha$ = $\frac{5}{15}$ in Equation (\ref{eq:param4}), we obtain the sixth-order compact derivative,denoted by \textbf{CD6} and written as:

\begin{equation}
    \frac{1}{3} \mathbf{ {U}}_{j-1}^{\prime}+\mathbf{ {U}}_{j}^{\prime}+\frac{1}{3} \mathbf{ {U}}_{j+1}^{\prime}=\frac{1}{36 \Delta x}\left(\mathbf{\hat {U}}_{j+2}-\mathbf{\hat {U}}_{j-2}\right)+\frac{7}{9 \Delta x}\left(\mathbf{\hat {U}}_{j+1}-\mathbf{\hat {U}}_{j-1}\right), \label{eq:cd6}
\end{equation}

\noindent where $j= 1, 2 , 3, .....,N-1$. Unlike the second-order central differences given by the Equations (\ref{eqn:expligrad}), which depend only on values at $j-1$, $j$ and $j+1$, compact finite differences depend on all the nodal values of the domain and therefore mimic the global dependence of the spectral methods. This global dependence results in a tridiagonal system of equations that the Thomas algorithm can easily invert. Near the boundary cells, lower-order one-sided difference formulas are used to approximate derivatives $\mathbf{{U}_0}'$ and $\mathbf{{U}_N}'$. The following third-order formulas are considered for both the CD4 and CD6 schemes in the present work.

\begin{equation}
\mathbf{{U}'_0}+2\mathbf{{U}'_1} = \frac{1}{\Delta x}(\frac{-5}{2}\mathbf{\hat {U}}_{0}+2\mathbf{\hat {U}}_1+\frac{1}{2}\mathbf{\hat {U}}_{2}),
\end{equation}  

\begin{equation}
\mathbf{{U}'_N}+2\mathbf{{U}'_{N-1}} = \frac{1}{\Delta x}(\frac{5}{2}\mathbf{\hat {U}}_{N}+2\mathbf{\hat {U}}_{N-1}+\frac{1}{2}\mathbf{\hat {U}}_{N-2}),
\end{equation}

{\color{black} The above formulas for the boundary cells is taken from \cite{moin2010fundamentals}, see their Equations 2.17.} For the computation of the second derivatives (Hessians), ${\mathbf{U}''_j}$, we compute the \textit{derivative of the first derivatives} obtained from Equations (\ref{eq:cd4}) and (\ref{eq:cd6}). For CD4:

\begin{equation}
    \frac{5}{14} \mathbf{ {U}''}_{j-1}+\mathbf{ {U}''}_{j}+\frac{5}{14} \mathbf{ {U}''}_{j+1}=\frac{b_1}{4 \Delta x}\left(\mathbf{ {U}'}_{j+2}-\mathbf{ {U}'}_{j-2}\right)+\frac{a_1}{2 \Delta x}\left(\mathbf{ {U}'}_{j+1}-\mathbf{ {U}'}_{j-1}\right), \label{eq:cd42}
\end{equation}

\noindent For CD6: 

\begin{equation}
    \frac{1}{3} \mathbf{ {U}''}_{j-1}+\mathbf{ {U}''}_{j}+\frac{1}{3} \mathbf{ {U}''}_{j+1}=\frac{1}{36 \Delta x}\left(\mathbf{ {U}'}_{j+2}-\mathbf{{U}'}_{j-2}\right)+\frac{7}{9 \Delta x}\left(\mathbf{ {U}'}_{j+1}-\mathbf{ {U}'}_{j-1}\right) \label{eq:cd62}
\end{equation}

{\color{black} It is emphasized that the orders of accuracy of these formulas for derivatives do not necessarily determine the order of accuracy of the resulting scheme because it is determined by the truncation error of the final residual, not by errors committed at an intermediate step (see a discussion on accuracy of fluxes and accuracy of the final residual in \textcolor{black}{references} \cite{VanLeerNishikawa_UltimateUnderstanding:JCP2021}). For the same reason, the use of the kappa scheme in the way we implement it does not imply third-order accuracy, and actually fourth-order accuracy is achieved for linear equations in the final residual as we will show later. }

After obtaining the derivatives ${ \mathbf{U}'}$ and ${ \mathbf{U}''}$, they are substituted into Equations (\ref{eqn:left_right}) to obtain the left and right reconstructed values necessary for the Riemann problem, completing the implicit gradient approach. Using CD4 or CD6 for the derivatives in Equation \ref{quadratic_polynomial_at_j} gives IG4 and IG6, respectively. In either case, the left- and right-interface values necessary for the desired approximate Riemann solver are, respectively:

\begin{equation}\label{eqn:IG}
\small
\begin{cases} 
 \mathbf{U}_{j+1 / 2}^{L,IG} &=\mathbf{\hat {U}}_{j}+\frac{1}{2} \mathbf{U}'_j+\frac{1}{12} \mathbf{U}''_j \\
\\											
\mathbf{U}_{j-1 / 2}^{R,IG} &=\mathbf{\hat {U}}_{j}-\frac{1}{2}  \mathbf{U}'_j+\frac{1}{12}  \mathbf{U}''_j 
\end{cases}
\begin{aligned}
\rightarrow &\text{$\textbf{U}'_j$ and $\textbf{U}''_j$ computed}\\& \text{by Eqns (\ref{eq:cd4}) and (\ref{eq:cd42}) is \textbf{IG4} and Eqns (\ref{eq:cd6}) and (\ref{eq:cd62}) is \textbf{IG6}}.
\end{aligned}
\end{equation}

 It is emphasized that both the kappa scheme and the implicit schemes for the derivatives are known formulas, but their combination had not been explored to the best of the authors' knowledge. In particular, as we will prove in the next section, it results in fourth-order upwind schemes {\color{black} with significantly low dispersive and dissipative errors.}

\subsubsection{\textcolor{black}{
{\color{black} Accuracy and properties of IG4 and IG6} }}
\label{sec-3.1.6}

{\color{black} Generally, schemes based on the quadratic reconstruction are third-order accurate at best. However, the use of high-order gradients makes it possible to go beyond third-order accuracy. To prove this, we apply the schemes to a linear convection equation $Q_t + F_x = 0$ with $F=Q$,
\begin{equation}\label{fourier:residual}
\frac{d \hat{Q}_j}{d t} = 
- \frac{1}{\Delta x} [ F_{j+1/2} - F_{j-1/2}  ],
\end{equation}
where
\begin{equation}\label{fourier:upwind_flux}
 F_{j+1/2} 
= 
{\color{black} 
\frac{1}{2} 
\left[  F( {Q}_{j+1/2}^R) + F({Q}_{j+1/2}^L)  \right]
-
 \frac{1}{2} ({Q}_{j+1/2}^R-{Q}_{j+1/2}^L)
 =
 {Q}_{j+1/2}^L,}
\end{equation}
and perform a Fourier analysis, which is especially useful for analyzing schemes with implicit gradients.   
Consider a Fourier mode: $\hat{Q}_\beta = \hat{Q}_0 \exp(i \beta  x / \Delta x)$, where $\hat{Q}_0$ is the amplitude, $\beta$ is the frequency, $\Delta x$ is the mesh spacing of a uniform grid, and $i=\sqrt{-1}$. Substituting it into the semi-discrete form (\ref{fourier:residual}), we find 
\begin{equation} 
\frac{d \hat{Q}_0}{d t} = {\cal F}^{exact}\hat{Q}_0  , 
\end{equation}
where ${\cal F}^{exact}$ denotes the exact convection operator:
\begin{equation}\label{fourier:exact_operator} 
{\cal F}^{exact}  = - \frac{ i \beta }{\Delta x}.
\end{equation}
This is the operator approximated by numerical schemes. Below, we derive the corresponding operators for the explicit and implicit gradients schemes, and measure the error by the leading deviation from the exact operator in the expansion for smooth components. The Fourier analysis is also useful for analyzing the dispersion and dissipative properties, which we will discuss subsequently.

\vspace{0.25cm}
\noindent{\it  Explicit gradients: Baseline scheme:}

First, we consider the explicit scheme with Equation (\ref{eqn:3linear}). Substituting the Fourier mode into the residual computed with the explicit gradients, we obtain
\begin{equation}\label{fourier:EG_semi_discrete}
\frac{d \hat{Q}_0}{d t} = {\cal F}^{EG}\hat{Q}_0  , 
\end{equation}
where
\begin{equation}\label{fourier:EG_res}
 {\cal F}^{EG}
 = 
- \frac{ ( \cos \beta -1 )^2 }{3} 
+ \frac{\sin \beta  ( \cos \beta - 4 )}{3}  i.
\end{equation} 
Expanding it for a small $\beta \approx O(\Delta x)$, we obtain
\begin{equation}\label{fourier:EG_res_expanded}
 {\cal F}^{EG}
 = 
- \frac{\beta }{\Delta x} 
\left[
i
+ \frac{\beta^3 }{12} 
+ \frac{ \beta^4}{30}i
- \frac{\beta^5}{72}
+ \frac{ \beta^6}{252}i
+ 
\cdots
\right],
\end{equation} 
which shows, compared with the exact operator (\ref{fourier:exact_operator}), that the leading error is third-order as expected: $O(\beta^3) = O(\Delta x^3)$. The third-order error appears in the real part, which indicates the error is dissipative. The leading dispersive error is fourth-order.

\vspace{0.25cm}
\noindent{\it Implicit gradients: IG4 and IG6:}

For the IG4 and IG6 schemes, we first derive the first and second derivatives. Consider a Fourier mode for the first derivative: $Q_x^{\beta} = G_0 \exp(i \beta/\Delta x)$, where $G_0$ is the amplitude, and substitute it into CD4 (Equation \ref{eq:cd4}) to obtain:
\begin{equation}\label{fourier:G_CD4}
\frac{5}{7} \cos \beta  G_0 + G_0 = 
\frac{ \sin \beta \cos \beta + 11 \sin \beta }
{7 \Delta x} Q_0 i , 
\end{equation} 
which can be solved for $G_0$: 
\begin{equation}\label{fourier:G_CD4_G0}
G_0 = 
\frac{  \sin \beta \cos \beta + 11  }{\Delta x (  5 \cos \beta + 7 ) i } Q_0 
.
\end{equation} 
Similarly, for the second derivative, substituting $Q_{xx}^{\beta} = H_0 \exp(i \beta/\Delta x)$, where $H_0$ is the amplitude, into CD4 (Equation \ref{eq:cd42}), solving it for $H_0$, and substituting Equation (\ref{fourier:G_CD4_G0}), we obtain
\begin{equation}\label{fourier:G_CD4_H0}
H_0 = - 
\frac{  (\sin \beta \cos \beta + 11 )^2 }{\Delta x^2 (  5 \cos \beta + 7 )^2} Q_0
.
\end{equation} 
Finally, substituting the Fourier modes, $Q^{\beta}$, $Q_{x}^{\beta}$, and $Q_{xx}^{\beta}$, into the IG4 scheme and eliminating $G_0$ and $H_0$ by the above equations, we obtain
\begin{equation}\label{fourier:IG4_res}
\begin{aligned}
 {\cal F}^{IG4} 
&=
\frac{( \cos \beta - 1 )^2 (\cos^3 \beta - 7 \cos^2 \beta + 11 \cos \beta - 5)}{12 (  5 \cos \beta + 7  )^2}
\\
&-
\frac{ \sin \beta  ( \cos^4 \beta  - 8 \cos^3 \beta + 78 \cos^2 \beta + 728 \cos \beta + 929)}{12 (  5 \cos \beta + 7  )^2} i,
\end{aligned}
\end{equation} 
which is expanded as  \begin{equation}\label{fourier:IG4_res_expanded}
 {\cal F}^{IG4} = 
- \frac{\beta }{\Delta x} 
\left[
i
+ \frac{\beta^4 }{720} i 
+ \frac{ \beta^6}{12096}i
+ \frac{\beta^7}{6912}
+ 
\cdots
\right].
\end{equation}
Clearly, it shows that it is fourth-order accurate with a leading fourth-order dispersive error and a sixth-order dissipative error. Note that there are no third- and fifth-order dissipative errors, i.e., no real part of $O(\beta^3)$ and $O(\beta^5)$, in contrast to the result for the explicit scheme  (\ref{fourier:EG_res_expanded}). It indicates that the use of implicit gradients effectively removed the third- and fifth-order dissipative errors and produces a significantly low-dissipation scheme.

Similarly, for the IG6 scheme, we first derive the amplitudes of the Fourier modes for the gradient and the second derivative by substituting the Fourier modes and solving the resulting equations for $G_0$ and $H_0$. The results are 
\begin{equation}\label{fourier:G_CD4_G6}
G_0 = 
\frac{  \sin \beta ( \cos \beta + 14)  }{3 \Delta x (  2 \cos \beta + 3 ) i } Q_0 ,
\end{equation}
and
\begin{equation}\label{fourier:G_CD6_H0}
H_0 = - 
\frac{  \sin^2 \beta ( \cos \beta + 14)^2  }{9 \Delta x^2 (  2 \cos \beta + 3 )^2} Q_0
.
\end{equation} 
Then, substituting the Fourier modes, $Q^{\beta}$, $Q_{x}^{\beta}$, and $Q_{xx}^{\beta}$, into the IG6 scheme and eliminating $G_0$ and $H_0$ by the above equations, we obtain
\begin{equation}\label{fourier:IG6_res}
\begin{aligned}
 {\cal F}^{IG6} &=
\frac{( \cos \beta - 1 )^2 ( \cos^3 \beta - 7 \cos^2 \beta + 26 \cos \beta - 20)}{ 108 (  2 \cos \beta + 3  )^2}
\\
&-
\frac{ \sin \beta  ( \cos^4 \beta  - 8 \cos^3 \beta + 105 \cos^2 \beta + 1070 \cos \beta + 1532)}{108 (  2 \cos \beta +37  )^2} i,
\end{aligned}
\end{equation} 
which is expanded as  \begin{equation}\label{fourier:IG6_res_expanded}
 {\cal F}^{IG6}  = 
- \frac{\beta }{\Delta x} 
\left[
i
+ \frac{\beta^4 }{720} i 
+ \frac{\beta^5}{1440}
- \frac{ \beta^6}{5040}i
+ \frac{17 \beta^7}{86400}
+ 
\cdots
\right].
\end{equation}
Therefore, the IG6 scheme is also fourth-order accurate with the leading fourth-order dispersive error. Note that there is a fifth-order dissipative error, which does not exist in IG4. Therefore, IG4 is expected to be less diffusive than IG6. Again, we emphasize that the order of accuracy of the scheme with implicit gradients is not necessarily determined by either the order of the reconstruction polynomial or the order of accuracy of the gradient algorithms. As we have just shown, the order of accuracy of the implicit-gradient-based schemes are not intuitive and needs to be analyzed correctly to reveal their accuracy and properties.

Note that the accuracy analysis is valid only for linear equations. For nonlinear equations, all the schemes are second-order accurate at best as long as the averaged flux is evaluated with reconstructed solutions: $[ F(Q^L_{j+1/2})+F(Q^R_{j+1/2})]/2$ \cite{VanLeerNishikawa_UltimateUnderstanding:JCP2021,FalseAccuracyUMUSCL:CiCP2021}. To preserve the high order of accuracy for nonlinear equations, the averaged flux needs to be evaluated by flux reconstruction $[ F^L_{j+1/2}+F^R_{j+1/2}]/2$, where $F^L_{j+1/2}$ and $F^R_{j+1/2}$ are computed by a direct flux reconstruction scheme, e.g., by applying the kappa scheme to the fluxes. \textcolor{black}{See  Nishikawa} \cite{Nishikawa_FSR:2020} for efficient flux reconstruction techniques based on solution derivatives (not using flux derivatives). Alternatively, one may reinterpret the numerical solution values as cell averages and apply a single-point high-order flux quadrature formula proposed in Refs.\cite{Buchmuller2014, tamaki2017efficient} or a multiple-point flux quadrature formula with a dimension-by-dimension solution reconstruction \cite{titarev2004finite}. However, it has been demonstrated that linearly high-order schemes can produce solutions with dramatically higher resolution than conventional second-order schemes for practical turbulent-flow simulations at a lower cost \cite{yang_harris:AIAAJ2016,HQYANG:AIAA2013-2021,yang_harris:CCP2018,GarciaBarakos:IJNMF2018}. The same has been demonstrated for flows with shock waves \textcolor {black} {by Zhang et al. in} \cite{zhang2011order}.

\vspace{0.25cm}
\noindent{\it  Dispersion and dissipation properties:}

\textcolor{black}{
Fig. \ref{fig_disp} shows the dispersion and dissipation properties of the IG4, IG6, and MUSCL schemes, {\color{black}  i.e., imaginary and real parts of  Equations (\ref{fourier:IG4_res}), (\ref{fourier:IG6_res}), and (\ref{fourier:EG_res}), respectively}. {\color{black} Figure \ref{fig:dispersion} shows that the dispersion property of IG4 is better than IG6. The dispersion property of the IG4 scheme is superior to that of the IG6 scheme.} Compared with the MUSCL scheme, which uses explicit gradients, it can be seen that explicit gradients are much more dispersive and dissipative than the implicit gradients considered here.}

\begin{figure}[H]
\centering
\subfigure[Dispersion]{\includegraphics[width=0.45\textwidth]{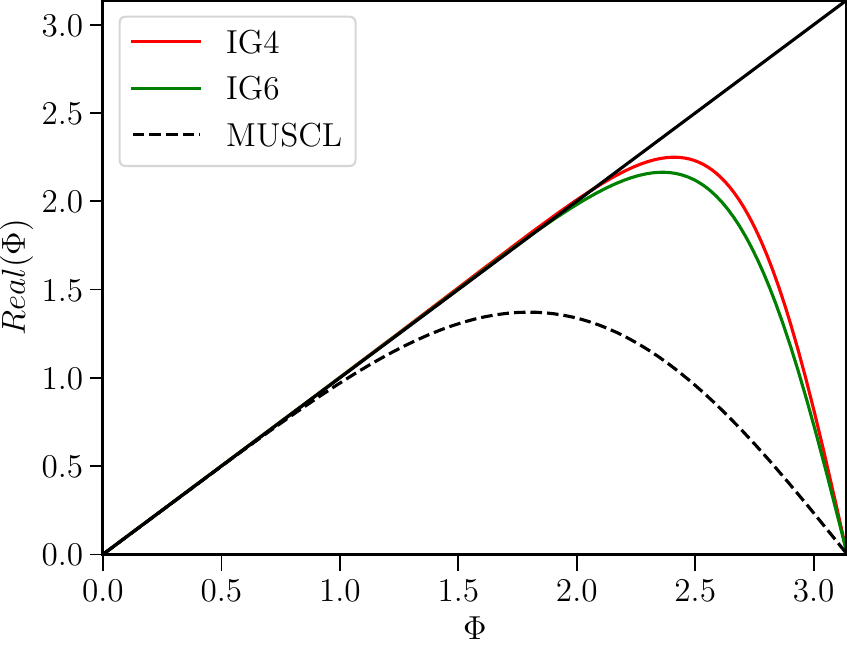}
\label{fig:dispersion}}
\subfigure[Dissipation]{\includegraphics[width=0.45\textwidth]{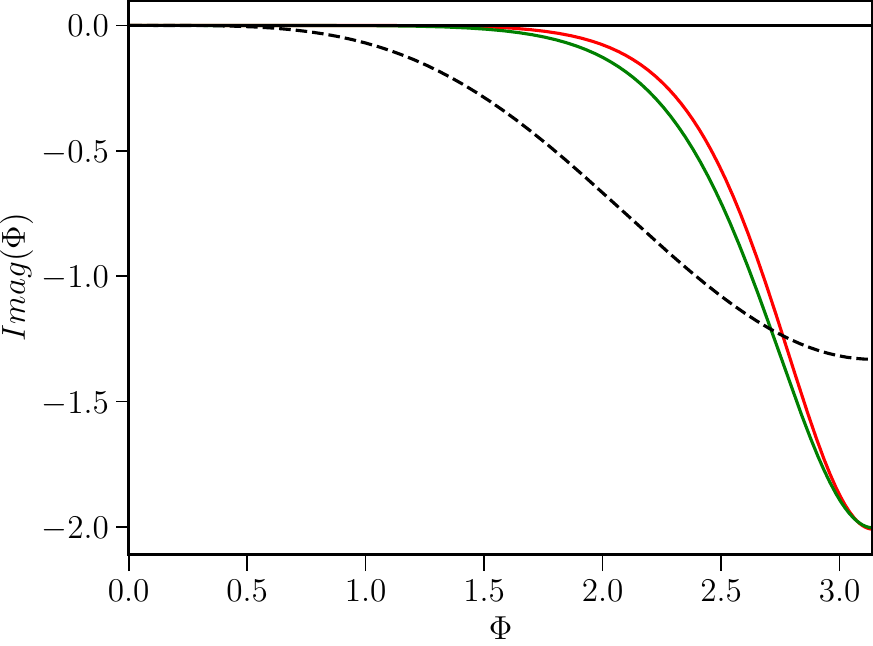}
\label{fig:dissipation}}
\caption{\textcolor{black}{Dispersion and Dissipation properties of the linear upwind schemes, where solid line: exact differentiation; green line: IG6; black dashed: MUSCL; red line: IG4.}}
\label{fig_disp}
\end{figure}
}

\begin{remark}\label{remark-recon}
\onehalfspacing
\normalfont {\textcolor{black}{In this paper, the proposed IG4 and IG6 schemes are termed gradient-based ``reconstruction'' schemes which is a misnomer. Instead, they should be called ``interpolation'' schemes as the ``reconstruction'' terminology is typically used for ``cell-averaged'' variables, and interpolation is used for ``point-values.'' Unfortunately, the first author has already used that terminology in \cite{CHAMARTHI2022105706}, which is already accepted. Therefore even in this paper, it is denoted as "reconstruction."}} 
\end{remark}

\subsection{Shock-capturing via BVD algorithm: {\color{black} IG4MP and IG6MP} }\label{sec-3.2}
The novel schemes, IG4 and IG6, derived in the earlier section, are linear in nature and therefore lead to oscillations for flows involving material interfaces and shocks. In this section, we describe the shock-capturing scheme using the BVD algorithm, previously presented in \cite{chamarthi2021high}, extended to the IG schemes for both single and multi-component flows. By comparing two different polynomials, the BVD algorithm selects the reconstruction polynomial with minimum numerical dissipation by evaluating the Total Boundary Variation (TBV) given by Equation (\ref{Eq:TBV}) for each cell for each \textit{primitive variable}:

\begin{equation}\label{Eq:TBV}
\textcolor{black}{{TBV}_{j}=\big|\mathbf{U}_{j-\frac{1}{2}}^{L}-\mathbf{U}_{j-\frac{1}{2}}^{R}\big|+\big|\mathbf{U}_{j+\frac{1}{2}}^{L}-\mathbf{U}_{j+\frac{1}{2}}^{R} \big|.}
\end{equation} 

For a given cell, the terms on the right-hand side of Equation (\ref{Eq:TBV}) represent the amount of numerical dissipation introduced in the numerical flux in Equation (\ref{eqn:Riemann}). The BVD algorithm compares the TBVs of the concerned polynomials and selects the one that is least dissipative at an interface. In the smooth regions of the flow, the IG6 or IG4 linear schemes will be used, and in the presence of discontinuities, the BVD algorithm will turn to the MP5 scheme \cite{suresh1997accurate}. The combination of IG4 and MP5 is denoted as \textbf{IG4MP} and, the combination of IG6 and MP5 is denoted as \textbf{IG6MP} in this paper, respectively. Both schemes are denoted together as \textbf{IGMP} schemes.}

In the following, the reconstruction procedure in the $x$-direction is discussed. Due to the dimension-by-dimension approach, the other directions are handled the same way and the procedure is summarized below:

\begin{description}
\item[Step i.] Evaluate the interface values by using the implicit gradient approach:
\begin{enumerate}[(a)]
\item {\color{black}Compute the first and second derivatives of the primitive variables {\color{black} by implicit gradient methods.} }
\item {\color{black}Form the reconstructed states by substituting the derivatives in Equation (\ref{eqn:IG})}, for each primitive variable (note, the superscript IG refers to either IG4 or IG6):

\begin{equation}\label{eqn:IG2}
\small
\begin{cases} 
 \mathbf{U}_{j+1 / 2}^{L,IG} &= \mathbf{\hat U}_{j}+\frac{1}{2} \mathbf{U}'_j+\frac{1}{12} \mathbf{U}''_j \\
\\											
\mathbf{U}_{j-1 / 2}^{R,IG} &= \mathbf{\hat U}_{j}-\frac{1}{2}  \mathbf{U}'_j+\frac{1}{12}  \mathbf{U}''_j 
\end{cases}
\begin{aligned}
\end{aligned}
\end{equation}
\end{enumerate}
\item [Step ii.] {\color{black}Evaluate the interface values using  the MP5 scheme. The steps involved are as presented in Appendix \ref{sec-appb}.}
\item[Step iii.] Calculate the TBV values for each cell $I_{j}$ for each candidate scheme,
\begin{equation}\label{Eq:TBVIG4}
\textcolor{black}{{TBV}_{j}^{IG}=\big|\mathbf{U}_{j-\frac{1}{2}}^{L,IG}-\mathbf{U}_{j-\frac{1}{2}}^{R,IG}\big|+\big|\mathbf{U}_{j+\frac{1}{2}}^{L,IG}-\mathbf{U}_{j+\frac{1}{2}}^{R,IG} \big|,}
\end{equation} 

\begin{equation}\label{Eq:TBVSC}
\textcolor{black}{{TBV}_{j}^{MP5}=\big|\mathbf{U}_{j-\frac{1}{2}}^{L,MP5}-\mathbf{U}_{j-\frac{1}{2}}^{R,MP5}\big|+\big|\mathbf{U}_{j+\frac{1}{2}}^{L,MP5}-\mathbf{U}_{j+\frac{1}{2}}^{R,MP5} \big|.}
\end{equation}

\item[Step iv.] \textcolor{black}{{\color{black}Modify} the interface values at $j-\frac{3}{2}$, $j-\frac{1}{2}$, $j+\frac{1}{2}$, and $j+\frac{3}{2}$} according to the following algorithm:

\begin{equation}
\text{if} \; {TBV}^{MP5} < {TBV}^{IG} \to \;  \left\{\begin{matrix}
\mathbf{U}^{K,IG}_{j-\frac{3}{2}} = \mathbf{U}^{{K},MP5}_{j-\frac{3}{2}}, \\ 
\\
\mathbf{U}^{K,IG}_{j-\frac{1}{2}} = \mathbf{U}^{{K},MP5}_{j-\frac{1}{2}}, \\ 
\\
\mathbf{U}^{K,IG}_{j+\frac{1}{2}} = \mathbf{U}^{{K},MP5}_{j-\frac{1}{2}}, \\ 
\\
\mathbf{U}^{K,IG}_{j+\frac{3}{2}} = \mathbf{U}^{{K},MP5}_{j+\frac{3}{2}},
\end{matrix}\right.
\end{equation}

where $K$ = $L$ or $R$.

\item[Step v.] \textcolor{black}{For each reconstructed state, $K$, we check the following conditions for positivity of pressure and density}

\begin{equation}
\begin{aligned}
\rho^{K} & \leq 0 \quad \text{\textbf{or}} \quad p^{K} \leq 0.
\end{aligned}
\end{equation}

 If the conditions given by Equations (\ref{eqn:single-positivity}) are not satisfied, the procedure is as follows,
\begin{equation}\label{eqn:single-positivity}
\mathbf{U}^K=\left\{\begin{array}{ll}
\mathbf{U}^{K, \text { IGMP}} & \text { default }, \\
\mathbf{U}^{K, \text { MP5}}& \text { if } \mathbf{U}^{K, \text { IGMP}} \textbf{ fails}, \\
\mathbf{U}^{K, \text { FO}} & \text { if } \mathbf{U}^{K, \text { MP5}} \textbf { fails },
\end{array}\right.
\end{equation}
where, $\mathbf{U}^{K, \text { FO}}$ is the first order approximation computed using Equation (\ref{eqn:first-order}).

\item[Step vi.] Evaluate the conservative variables, $(\mathbf{Q}_{j+\frac{1}{2}}^{L},\mathbf{Q}_{j+\frac{1}{2}}^{R})$, from the primitive variables, $(\mathbf{U}_{j+\frac{1}{2}}^{L},\mathbf{U}_{j+\frac{1}{2}}^{R})$ obtained from the above procedure, and compute the interface flux $\mathbf{F}^{\rm Riemann}_{j+\frac{1}{2}}$
\end{description}

\begin{remark}\label{remark-3.5}
\onehalfspacing
\normalfont The choice of MP5 scheme as a candidate polynomial for shock-capturing has been arrived at by testing several different combinations of linear and nonlinear schemes. It is also possible to use a third-order reconstruction with the minmod limiter and THINC, and readers are referred to \textcolor{black}{Appendix A and B of \cite{chamarthi2021high}. The BVD algorithm used here is similar to the one used by that of Chamarthi and Frankel \cite{chamarthi2021high}.} 
\end{remark}

\begin{remark}\label{remark-posi}
\onehalfspacing
\normalfont \textcolor{black}{In this paper, the positivity preserving approach given by Equation (\ref{eqn:single-positivity}) was activated only for the test case in Example \ref{ex:dmr}. One may use even more sophisticated positivity preserving approaches available in the literature \cite{zhang2012positivity,hu2013positivity}.} 
\end{remark}

\subsection{Riemann solver}\label{sec-3.3}

Approximate Riemann solvers approximate the convective flux after obtaining reconstructed states at the interface as explained in the earlier section. This section illustrates how the HLLC approximate Riemann solver \cite{batten1997choice,toro1994restoration} approximates convective fluxes. For simplicity, only the HLLC approximations for both single and multicomponent flows in the x-direction are illustrated in this section. The HLLC flux in $x$-direction is given by:

\begin{equation}\label{eq:HLLC}
\mathbf{F}^{\rm Riemann}= \mathbf{F}^{HLLC}=\left\{\begin{array}{ll}
\mathbf{F}_{L} & , \text { if } \quad 0 \leq S_{L}, \\
\mathbf{F}_{* L} & , \text { if } \quad S_{L} \leq 0 \leq S_{*}, \\
\mathbf{F}_{* R} & , \text { if } S_{*} \leq 0 \leq S_{R}, \\
\mathbf{F}_{R} & , \text { if } \quad 0 \geq S_{R},
\end{array}\right.
\end{equation}
\begin{equation}
\mathbf{F}_{* K}=\mathbf{F}_{K}+S_{K}\left(\mathbf{Q}_{* K}-\mathbf{Q}_{K}\right),
\end{equation}
where $L$ and $R$ are the left and right states respectively. With $K$ = $L$ or $R$, the star state \textcolor{black}{quantities are} defined as:

\begin{equation}\label{eqn-hllc-single}
\mathbf{Q}_{* K}=\left(\frac{S_{\mathrm{K}}-u_{\mathrm{K}}}{S_{\mathrm{K}}-S_{*}}\right)\left[\begin{array}{c}
\rho_{K} \\
\rho_{\mathrm{K}} S_{*} \\
\rho_{K} v_{K} \\
\rho_{K} w_{K} \\
E_{k}+\left(S_{*}-u_{K}\right)\left(\rho_{K} S_{*}+\frac{p_{K}}{S_{K}-u_{K}}\right)
\end{array}\right].
\end{equation}

 In the above expressions, the waves speeds  $S_L$ and $S_R$ can be obtained as suggested by Einfeldt \cite{einfeldt1988godunov}, $ S_L = min(u_{L}-c_L, \tilde{u}-\tilde{c}) \ \text{and} \ S_R = max(u_{R}+c_R, \tilde{u} +\tilde{c})$, where $\tilde {u}$ and $\tilde {c}$ are the Roe averages from the left and right states (for the definitions of the Roe averaged quantities, see Blazek \cite{blazek2015computational}). Batten et al. \cite{batten1997choice} provided a closed form expression for $S_*$:

\begin{align}
 \label{eqn:HLLCmiddlewaveestimate}
 S_* = \frac{p_R - p_L + \rho_Lu_{nL}(S_L - u_{nL}) - \rho_Ru_{nR}(S_R - u_{nR})}{\rho_L(S_L - u_{nL}) -\rho_R(S_R - u_{nR})}.
\end{align}

\subsection{Temporal integration}\label{sec-3.4}
Finally, the semi-discrete approximation of the governing equations is temporally integrated. The conserved variables are  integrated in time using the following third-order total-variation-diminishing (TVD) Runge-Kutta scheme \cite{jiang1995}:

\begin{eqnarray} \label{rk}
\mathbf{\hat Q}^{(1)} & = & \mathbf{\hat Q}^n + \Delta t \mathbf{Res}(\mathbf{\hat Q}^n), \nonumber \\
\mathbf{\hat Q}^{(2)} & = & \frac{3}{4} \mathbf{\hat Q}^n + \frac{1}{4} \mathbf{\hat Q}^{(1)} + \frac{1}{4}\Delta t \mathbf{Res}(\mathbf{\hat Q})^{(1)} ,\\
\mathbf{\hat Q}^{n+1} & = & \frac{1}{3} \mathbf{\hat Q}^n + \frac{2}{3} \mathbf{\hat Q}^{(2)} + \frac{2}{3}\Delta t \mathbf{Res}(\mathbf{\hat Q}^{(2)}), \nonumber
\end{eqnarray}
where {\color{black} the $(j,i,k)$ component of $\mathbf{Res}$ is given by the right-hand side of Equation (\ref{eqn-differencing}). }
The superscripts ${n}$ and ${n+1}$ denote the current and the subsequent time-steps, and superscripts ${(1)-(2)}$ correspond to intermediate steps. The time step $\Delta t$ is computed as:

\begin{equation}\label{eqn:cfl}
\Delta t= \text{CFL} \cdot \min \left(\Delta t_{viscous}, \Delta t_{inviscid}\right),
\end{equation}
where
\begin{equation}
\Delta t_{inviscid}=  \min _{j,i,k}\left(\frac{\Delta x_{i}}{\left|u_{j,i,k}\right|+c_{j,i,k}}, \frac{\Delta y_{j}}{\left|v_{j,i,k}\right|+c_{j,i,k}}, \frac{\Delta z_{k}}{\left|w_{j,i,k}\right|+c_{j,i,k}}\right),
\end{equation}

\noindent where $c$ is the speed of sound and given by $c=\sqrt{\gamma{p}/\rho}$, and 
\textcolor{black}{\begin{equation}
 \Delta t_{viscous}=  \min _{j,i,k}\left(\frac{1}{\alpha} \frac{\Delta x_{i}^{2}}{\nu_{j,i,k}}, \frac{1}{\alpha}  \frac{\Delta y_{j}^{2}}{\nu_{j,i,k}}, \frac{1}{\alpha}  \frac{\Delta z_{k}^{2}}{\nu_{j,i,k}}\right),
\end{equation}}

\noindent where $\alpha = 4$ (see Chamarthi et al. \cite{chamarthi2022} for the chosen constant of $\alpha$) and $\nu = \mu / \rho$ is the kinematic viscosity. The details of the time-step restriction are presented in \cite{chamextend,chamarthi2022}. Time integration is performed with a CFL = 0.2 for all the test cases.

\section{Results}\label{sec-4}

In this section, we present various results displaying the performance of the proposed schemes. We compare the IGMP methods with the TENO5 scheme of Fu et al. \cite{fu2019low}. In the following, the IGMP schemes use primitive variable reconstruction, for the reasons delineated above, while the TENO5 scheme uses conservative variable reconstruction as the cited reference presents.

\subsection{Order of accuracy}
\begin{example}\label{euler-accuracy}{Accuracy of the proposed schemes}
\end{example}
\textcolor{black}{Firstly}, we present the order of accuracy (OOA) of the new schemes by convecting the initial profile given by equation \eqref{accu-euler} in the domain $x,y \in [-1, 1]$. The solution is obtained at time $t=2$, and the time-step is varied as a function of the grid size as \textcolor{black}{$\Delta t = \text{CFL} \Delta x$.}

\begin{align}\label{accu-euler}
(\rho,u,v,p)= (1+0.5 \sin(x+y),\ \ 1.0,\ \, 1.0 \ \ 1.0).
\end{align}
{\color{black} This exact solution effectively linearizes the Euler equations \cite{VanLeerNishikawa_UltimateUnderstanding:JCP2021,FalseAccuracyUMUSCL:CiCP2021} and therefore we expect fourth-order accuracy for the new schemes, which are linearly fourth-order accurate as shown above.} The $L_2$ norm of the error between the exact {\color{black} solution evaluated at the cell center} and the obtained solution is used to compute the OOA. Table \ref{tab:ooa} presents the OOAs obtained for this test case.

\begin{table}[H]
  \centering
   \caption{$L_2$ errors and numerical order of accuracy for the linear test case. $N$ is the number of cells in the domain}
\footnotesize
    \begin{tabular}{crccccc}
\hline
    N     & \multicolumn{1}{c}{TENO5} & OOA   & IG6MP & OOA   & IG4MP & OOA \\
\hline
    $10^2$ & 6.79E-03 & -     & 5.98E-04 & -     & 4.65E-04 & - \\
\hline
    $20^2$ & 2.24E-04 & 4.92  & 4.59E-05 & 3.71  & 4.37E-05 & 3.41 \\
\hline
    $40^2$ & 7.06E-06 & 4.98  & 2.54E-06 & 4.18  & 2.30E-06 & 4.25 \\
\hline
    $80^2$ & 2.21E-07 & 5.00  & 1.77E-07 & 3.84  & 1.74E-07 & 3.72 \\
\hline
    \end{tabular}%
 \label{tab:ooa}%
\end{table}%

\textcolor{black}{We have the following observations regarding order of accuracy:}

\begin{itemize}
\item \textcolor{black}{The order accuracy of the IG schemes is consistent with the mathematical proof provided in section \ref{sec-3.1.5}. Therefore it is possible to obtain fourth-order accuracy with the kappa scheme given by Equation (\ref{quadratic_polynomial_at_j}).}
	\item \textcolor{black}{$L_2$ norm errors given in Table \ref{tab:ooa} also indicate that the absolute error of the IGMP schemes is nearly the same as that of the fifth-order TENO5 scheme. However, \textcolor{black}{as it will be shown in the following examples}, the IGMP schemes will give better resolution than the TENO5 scheme. Hu et al. \cite{hu2012dispersion} optimized the linear schemes for favourable spectral properties satisfying the dispersion-dissipation relation but such optimization lead to order degeneration. Even though the proposed schemes are only fourth-order accurate, they have superior dispersion and dissipation properties.} 
\end{itemize}

\textcolor{black}{Next, we consider the nonlinear test case proposed by Yee et al.\cite{yee1999low}, where an isentropic vortex is convected in an inviscid free stream. The computations are carried out on a computational domain of [0, 10] $\times$ [0, 10] with periodic boundary conditions on all sides. The case is run until $t$=10. The initial conditions for this test case are:}

\begin{equation}
\begin{array}{l}
\rho=\left[1-\frac{(\gamma-1) \beta^{2}}{8 \gamma \pi^{2}} e^{\left(1-r^{2}\right)}\right]^{\frac{1}{p-1}}, \quad r^{2}=\bar{x}^{2}+\bar{y}^{2}, \\
(u, v)=(1,1)+\frac{\beta}{2 \pi} e^{\frac{1}{2}\left(1-r^{2}\right)}(-\bar{y}, \bar{x}), \quad \bar{x}=x-x_{v c}, \quad \bar{y}=y-y_{v c}, \\
p=\rho^{\gamma},
\end{array}
\end{equation}
\textcolor{black}{where $(x_{v c},y_{v c})$ = (5,5) are the coordinates of the center of the initial vortex and $\beta$ = 5. {\color{black} Note that this problem can still be linearized if the parameter $\beta$ is small, and the choice $\beta=5$ is made to avoid it.}}  See \textcolor{black}{references} \cite{VanLeerNishikawa_UltimateUnderstanding:JCP2021, FalseAccuracyUMUSCL:CiCP2021} for details. {\color{black} Results are summarized in Table  \ref{tab:addlabel}. We can observe clearly from Table \ref{tab:addlabel} that the proposed schemes are only second-order accurate for nonlinear cases, as expected. The TENO5 scheme, using conservative variables, is also second-order accurate as mentioned in \cite{fu2019low}.}

\begin{table}[H]
 \centering
  \footnotesize
 \caption{$L_2$ errors and numerical order of accuracy for the nonlinear test case. $N$ is the number of cells in the domain}
    \begin{tabular}{cccrcrc}
      \hline
    N     & TENO5 & OOA   & \multicolumn{1}{c}{IG6MP} & OOA   & \multicolumn{1}{c}{IG4MP} & OOA \\
          \hline
    $25^2$ & 4.33E-03 & -     & 3.14E-03 & -     & 3.14E-03 & - \\
           \hline
    $50^2$ & 3.21E-04 & 2.49  & 6.54E-04 & 2.26  & 6.55E-04 & 2.26 \\
	\hline
    $100^2$ & 3.30E-05 & 2.17  & 1.64E-04 & 2.00  & 1.64E-04 & 2.00 \\
	\hline
    $200^2$ & 8.35E-06 & 1.99  & 4.10E-05 & 2.00  & 4.10E-05 & 2.00 \\
	 \hline
    \end{tabular}%
  \label{tab:addlabel}%
\end{table}%

{\color{black} 
}

\textcolor{black}{Zhang et al. \cite{zhang2011order} demonstrated}  that the finite-volume WENO method with mid-point rule is only second-order accurate for nonlinear systems, and the Gaussian integral rule is necessary for high-order accuracy. However, they also noted that the resolution characteristics are often comparable for flows involving discontinuities despite the difference in the order of accuracy. In the present paper, we are also interested in flows involving discontinuities. {\color{black}
Also, such linearly high-order schemes have been demonstrated to serve as very low-dissipation schemes for practical turbulent-flow simulations 
 \cite{yang_harris:AIAAJ2016,HQYANG:AIAA2013-2021,yang_harris:CCP2018,GarciaBarakos:IJNMF2018}. In the following sections, we will demonstrate that the developed schemes do indeed serve as very  low-dissipation/dispersion schemes for strongly nonlinear problems with shock waves. 
 }

\subsection{One-dimensional Euler equations}

In this subsection, we consider the test cases for the one-dimensional Euler equations.

\begin{example}\label{sod}{Shock tube problems}
\end{example}

The two shock-tube problems proposed by Sod \cite{sod1978survey} and Lax \cite{lax1954weak}, are solved by the proposed schemes, IG6MP and IG4MP. The solutions are obtained by setting the specific heat ratio to be $\gamma=1.4$, and are compared with that of results from an exact Riemann solver \cite{toro2009riemann}. The initial conditions for the Sod test case and the Lax problem are given by the following initial conditions \eqref{sod_prob}, and \eqref{lax_prob}, respectively.

\begin{align}\label{sod_prob}
(\rho,u,p)=
\begin{cases}
(0.125,\ \ 0,\ \ 0.1),&\quad 0<x<0.5,\\
(1,\ \ 0,\ \ 1),&\quad 0.5 \leq x<1,
\end{cases}
\end{align}
\begin{align}\label{lax_prob}
(\rho,u,p)=
\begin{cases}
(0.445,\ \ 0.698,\ \ 3.528),&\quad 0<x<0.5,\\
(0.5,\ \ 0,\ \ 0.571),&\quad 0.5 \leq x<1.
\end{cases}
\end{align}

\begin{figure}[H]
\centering
\subfigure[Density]{\includegraphics[width=0.4\textwidth]{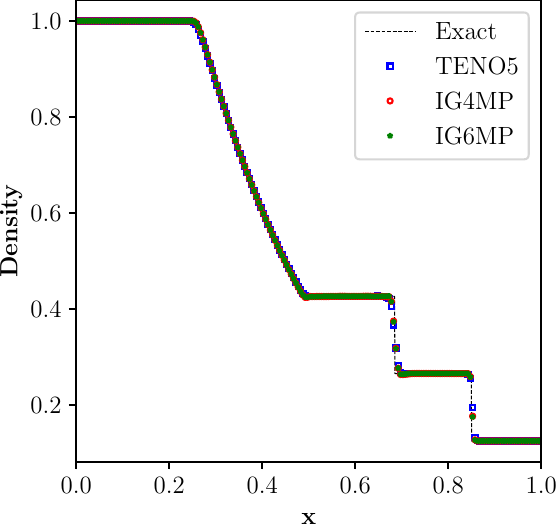}
\label{fig:sod-den}}
\subfigure[Pressure]{\includegraphics[width=0.4\textwidth]{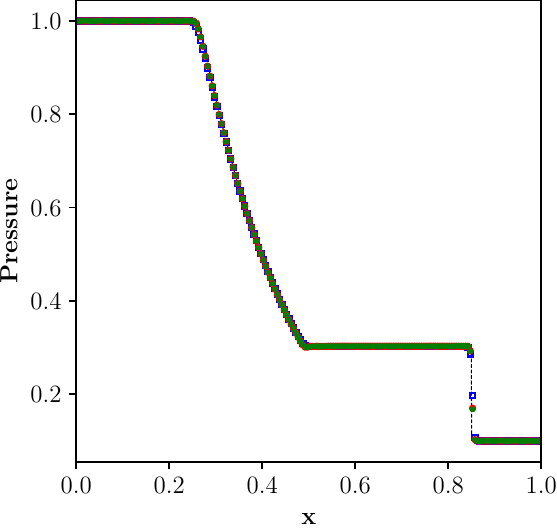}
\label{fig:sod-pres}}
\caption{Numerical solution for Sod problem in Example \ref{sod} for $N = 200$ points, where dashed line: reference solution; green stars: IG6MP; blue squares: TENO5; red circles: IG4MP. }
\label{fig_sod}
\end{figure}

First, the Sod test case is used to assess the shock-capturing ability of the scheme. The case is run using 200 cells until a final time of $t = 0.2$. Fig. \ref{fig_sod} shows the density and velocity profiles of the proposed schemes and the TENO5 scheme compared with the exact solution. The solutions by all the schemes are in good agreement with the exact solution in addition to the absence of overshoots in regions of discontinuities.

\begin{figure}[H]
\centering
\subfigure[Density]{\includegraphics[width=0.4\textwidth]{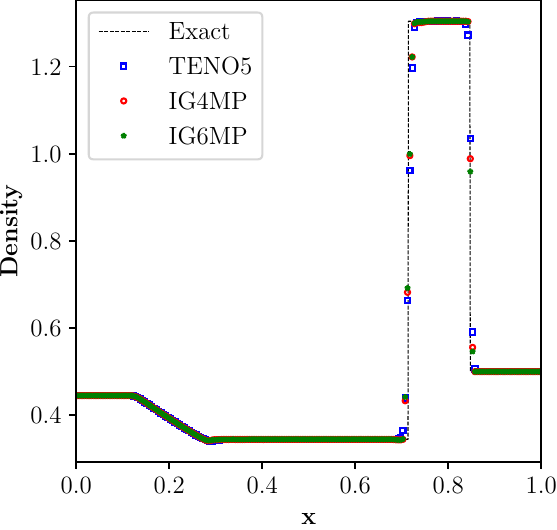}
\label{fig:lax-den}}
\subfigure[Velocity]{\includegraphics[width=0.4\textwidth]{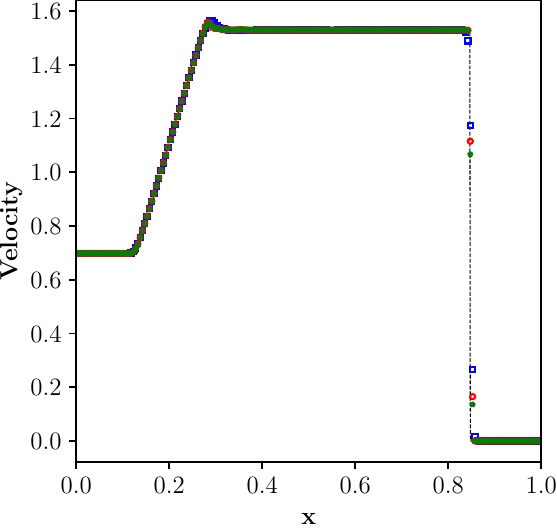}
\label{fig:lax-pres}}
\caption{Numerical solution for Lax problem in Example \ref{sod} for $N=200$ cells, where dashed line: reference solution; green stars: IG6MP; blue squares: TENO5; red circles: IG4MP.}
\label{fig_lax}
\end{figure}

Second, we used $200$ cells for the Lax problem, and the solution is obtained at time $t=0.14$. Fig. \ref{fig_lax} presents the density and velocity of the new schemes and the \textcolor{black}{TENO5} scheme compared with the exact solution. The IG4MP scheme resolves the features of the flow while avoiding oscillations. Observing the velocity profile, the proposed schemes can capture the discontinuity with relatively fewer cells compared to the TENO5 scheme.

\begin{example}\label{Shu-Osher}{Shu-Osher problem}
\end{example}

Third, we consider the Shu-Osher problem \cite{Shu1988}, which is a one-dimensional idealization of shock-turbulence interaction that simulates the interaction of a right moving shock wave for a given Mach number ($M=3$) superimposed with a perturbed density field. The initial conditions for this problem are:

\begin{align}\label{shock_den_prob}
(\rho,u,p)=
\begin{cases}
(3.857143,\ \ 2.629369,\ \ 10.3333),&\quad -5<x<-4,\\
(1+0.2\sin(5x),\ \ 0,\ \ 1),&\quad -4 \leq x<5.
\end{cases}
\end{align}

\begin{figure}[H]
\centering
\subfigure[Global profile]{\includegraphics[width=0.4\textwidth]{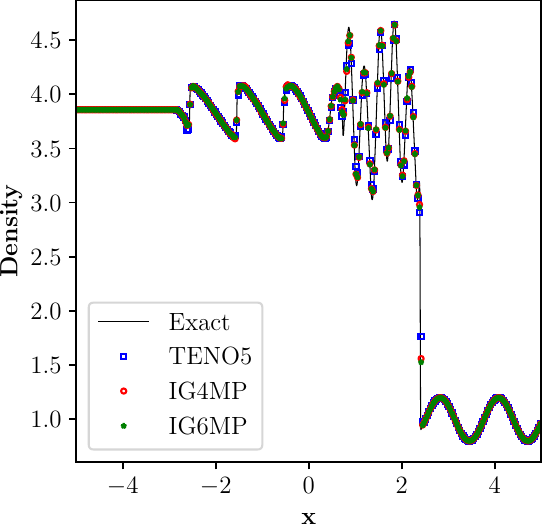}
\label{fig:shu-global}}
\subfigure[Local profile]{\includegraphics[width=0.42\textwidth]{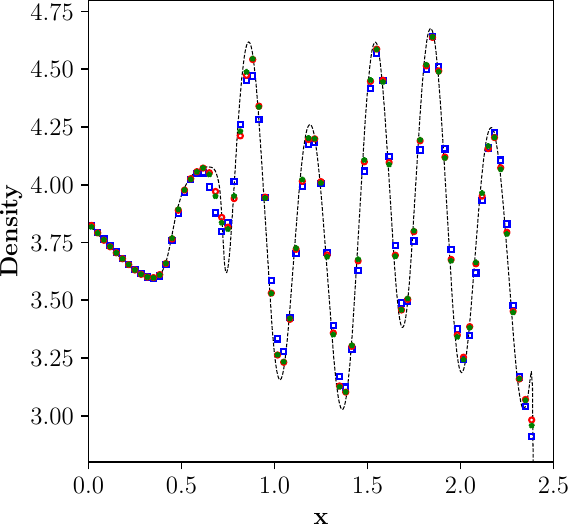}
\label{fig:shu-local}}
\caption{\textcolor{black}{Density profile for Shu-Osher problem, Example \ref{Shu-Osher}, \textcolor{black}{on a grid of $300$ cells}, where solid or dashed line: reference solution; green stars: IG6MP; blue squares: TENO5; red circles: IG4MP.}}
\label{fig:1d-SO}
\end{figure}

The solution is obtained for a time $t=1.8$ on a grid size of $300$ cells. The reference solution is obtained using the WENOZ scheme \cite{Borges2008} on a fine grid of $1600$ cells. The density profiles are shown in Fig. \ref{fig:1d-SO}. It is observed that both IGMP schemes perform well in capturing the post-shock oscillations in density. \textcolor{black}{Notably, both the proposed schemes capture the peaks and troughs of the density very well and better than the TENO5 scheme.}

\begin{example}\label{Titarev-Toro}{Titarev-Toro problem}
\end{example}

The last one-dimensional test case we consider is the shock-entropy wave problem of Titarev-Toro \cite{titarev2004finite}. In this test, a high-frequency sinusoidal wave interacts with a shock wave. The test case reflects the ability of the scheme to capture the extremely high-frequency waves. The initial conditions are given by equation \eqref{shock_tita} on a domain of $[-5,5]$,

\begin{align}\label{shock_tita}
(\rho,u,p)=
\begin{cases}
(1.515695,\ \ 0.523326,\ \ 1.805),&\quad x<-4.5,\\
(1+0.1\sin(20x\pi),\ \ 0,\ \ 1),&\quad x \geq -4.5.
\end{cases}
\end{align}

\begin{figure}[H]
\centering
\subfigure[Global profile]{\includegraphics[width=0.4\textwidth]{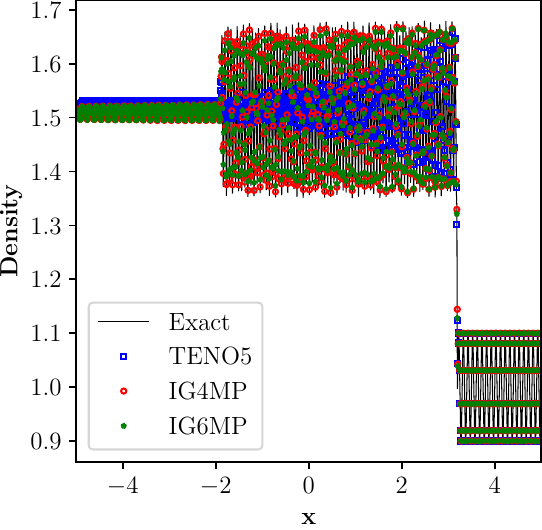}
\label{fig:tita1}}
\subfigure[Local profile]{\includegraphics[width=0.4\textwidth]{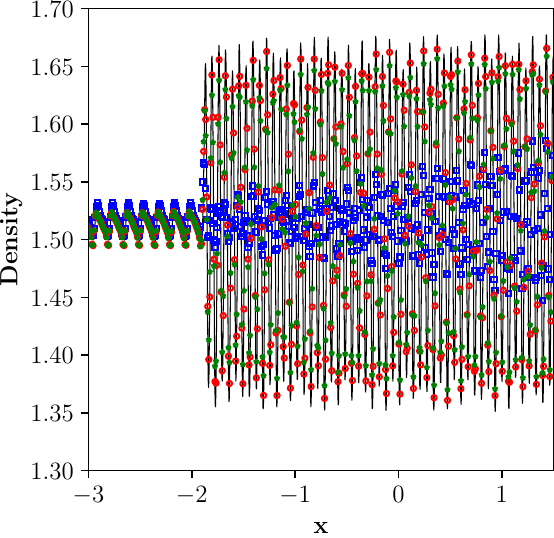}
\label{fig:tita2}}
\caption{Density profiles obtained by various schemes for Example \ref{Titarev-Toro}. Solutions are obtained for grid size of $1000$ cells. Solid line: reference solution; green stars: IG6MP; blue squares: TENO5; red circles: IG4MP.}
\label{fig_tita}
\end{figure}

 We used $1000$ cells to simulate the problem until time $t=5$ and compared with the reference solution obtained from a WENOZ simulation on a fine grid size of $3000$ cells. The results obtained presented in Fig. \ref{fig_tita} indicate that the IGMP schemes can accurately capture the high-frequency wave. Specifically, we observe that the IGMP schemes capture the linear region significantly better than the \textcolor{black}{TENO5} scheme.

\subsection{Multi-dimensional test cases for Euler equations}

\begin{figure}[H]
\centering\offinterlineskip
\subfigure[TENO5]{\includegraphics[width=0.48\textwidth]{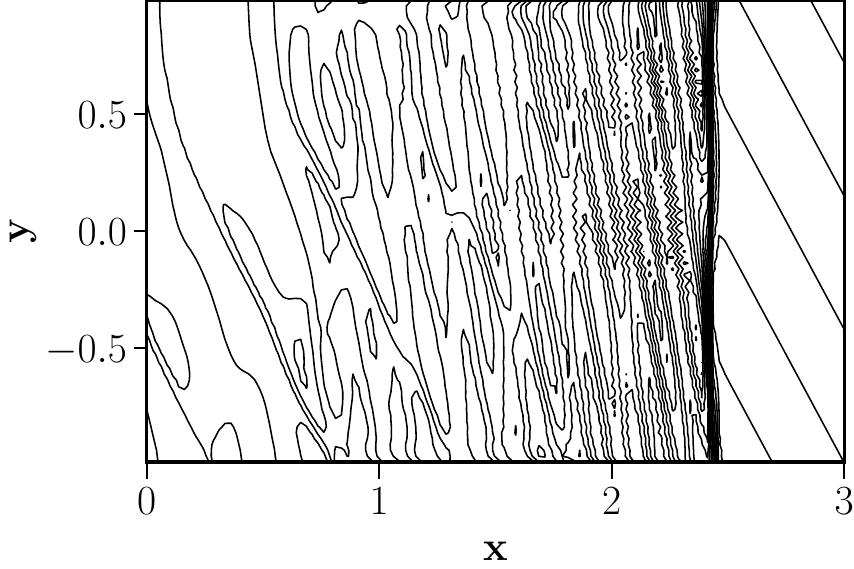}
\label{fig:SE-TENO5}}
\subfigure[IG6MP]{\includegraphics[width=0.48\textwidth]{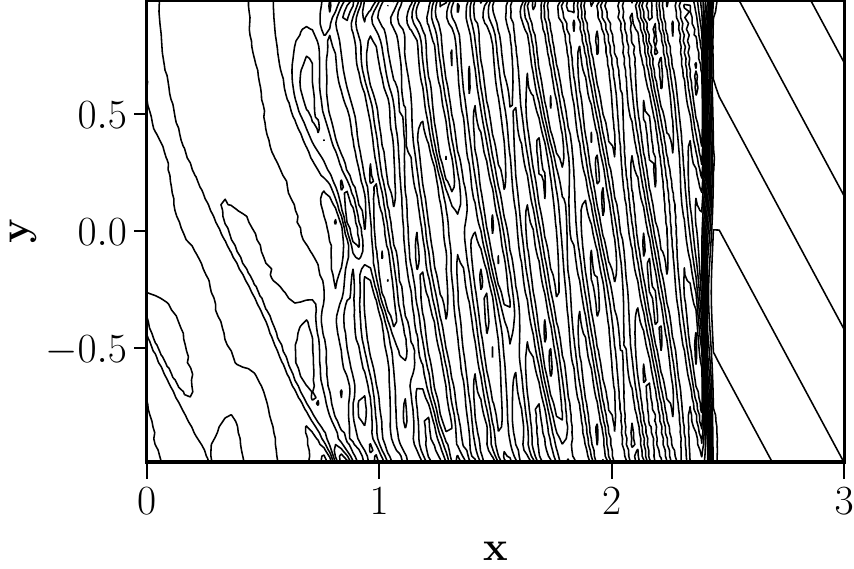}
\label{fig:SE-BVD5}}
\subfigure[IG4MP]{\includegraphics[width=0.48\textwidth]{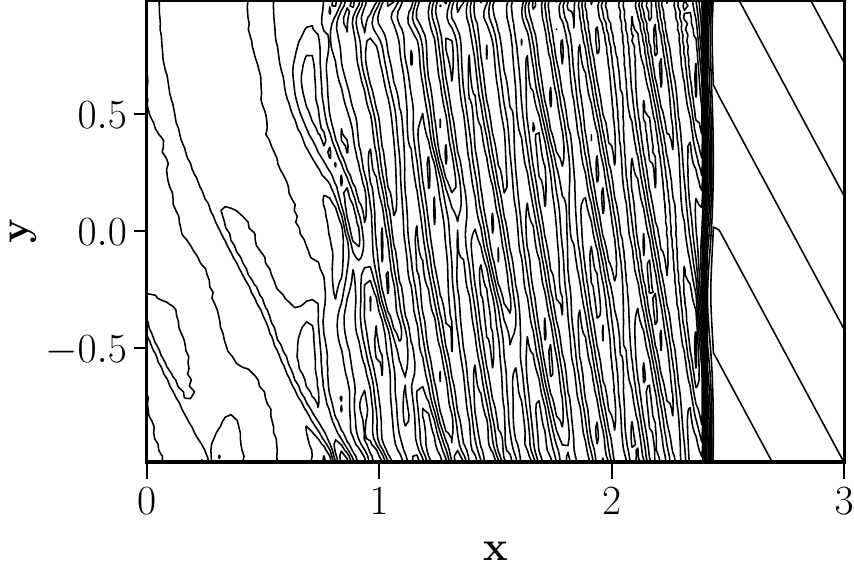}
\label{fig:SE-BVD}}
\subfigure[Local profile]{\includegraphics[width=0.48\textwidth]{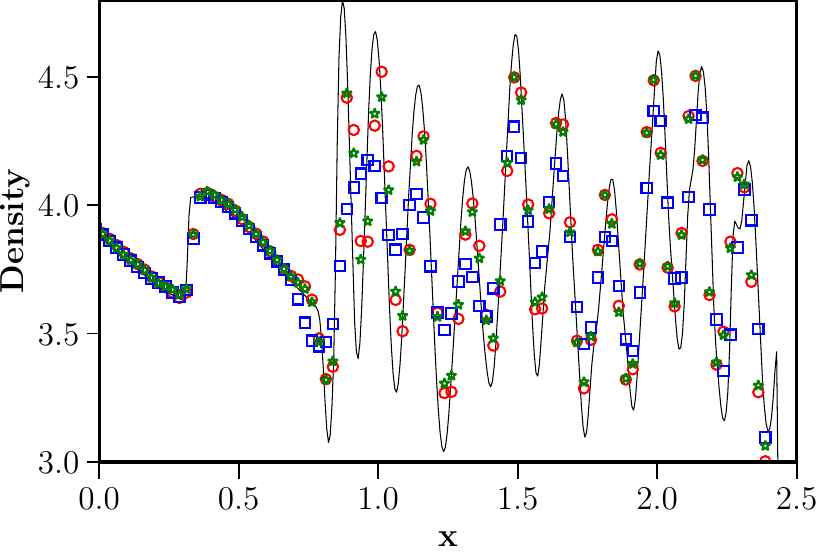}
\label{fig:SSE-Compare}}
\caption{Density contours for the 2D shock-entropy wave test at $t=1.8$, Example \ref{shock-entropy}, for various schemes are shown in Figs. (a), (b) and (c). Fig. (d) shows the local density profile in the region with high-frequency waves for all the schemes. Solid line: reference solution; green stars: IG6MP; blue squares: TENO5; red circles: IG4MP.}
\label{fig_shock_entropy}
\end{figure}

\begin{example}\label{shock-entropy}{2D Shock-Entropy Wave Test}
\end{example}
In this test case we consider the two-dimensional shock-entropy wave interaction problem proposed in \cite{acker2016improved}. The initial conditions for the test case are as follows,
\begin{align}\label{shock_entropy}
(\rho,u,v,p)=
\begin{cases}
(3.857143, \ \ 2.629369,\ \ 0,\ \ 10.3333),&\quad x\leq-4,\\
(1+0.2\sin(10x \cos\theta+10y\sin\theta),\ \ 0,\ \ 0,\ \ 1),&\quad otherwise,
\end{cases}
\end{align}
with $\theta$ = $\pi/6$ over a domain of $[-5,5]\times [-1,1]$. The initial sine waves make an angle of $\theta$ radians with the $x$ axis. Initial conditions are modified as in \cite{deng2019fifth} with a higher frequency for the initial sine waves compared to that of \cite{acker2016improved} to show the benefits of the proposed method. A grid of $400\times 80$ is chosen and the case is run until $t=1.8$. The reference solution is computed on a fine grid of $1600 \times 320$ using the WENOZ scheme. Density contour plots shown in Fig. \ref{fig_shock_entropy} indicate that the proposed schemes significantly improve the resolution of the flow structures. The local density profile along $y=0$ is presented in Fig. \ref{fig:SSE-Compare}. The results demonstrate that IG6MP and IG4MP retain the desirable shock-capturing features in the MP5 scheme while capturing the high-frequency region better than the \textcolor{black}{TENO5} scheme. Also, the TENO5 scheme has oscillations, shown in Fig. \ref{fig:SE-TENO5}, for this test case.

\begin{example}\label{ex:rp}{Riemann Problem }
\end{example}

In this test case we consider the Riemann problem of \cite{schulz1993numerical} described as configuration 3. The initial conditions of the problem are given by equation \eqref{riemann_problem} with constant states of the primitive variables along the lines $x=0.8$, and $y=0.8$ in the domain $x,y \in [0,1]$. This produces four shocks at the interfaces of the four quadrants. Also, the small-scale complex structures generated along the slip-lines due to the Kevin-Helmholtz instabilities serve to assess the numerical dissipation of the scheme. Non-reflective boundary conditions are employed on all four boundaries. The test case is run until $t=0.8$ on a grid of size $400 \times 400 $. 

 \begin{equation}\label{riemann_problem}
(\rho, u,v, p)=\left\{\begin{array}{ll}
        (1.5, 0, 0, 1.5), ~~~~~~~~~~~~~~~~~~~~~~~~~~~~~~~~\mbox{if} ~~x > 0.8, ~~y > 0.8, \\
        (33/62, 4/\sqrt{11}, 0, 0.3), ~~~~~~~~~~~~~~~~~~~~~\mbox{if} ~~ x\leq 0.8, ~~y > 0.8, \\
        (77/558, 4/\sqrt{11}, 4/\sqrt{11}, 9/310), ~~~~~~~~~\mbox{if} ~~ x \leq 0.8, ~~y\leq 0.8, \\
        (33/62, 0, 4/\sqrt{11}, 0.3), ~~~~~~~~~~~~~~~~~~~~~\mbox{if} ~~x > 0.8, ~~ y\leq 0.8.
        \end{array}\right.
\end{equation}

\begin{figure}[H]
\begin{onehalfspacing}
\centering\offinterlineskip
\subfigure[TENO5]{\includegraphics[width=0.48\textwidth]{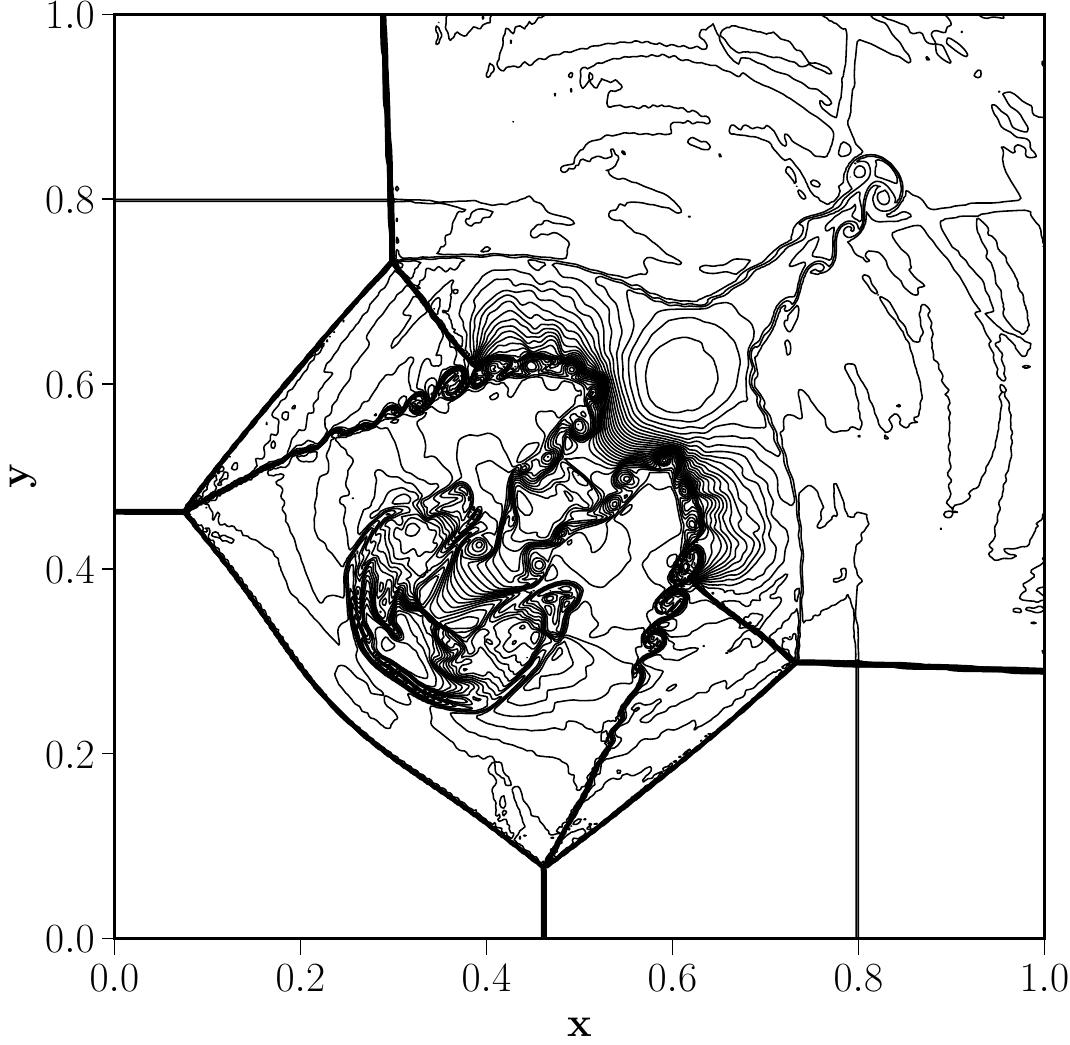}
\label{fig:MP5_RM}}
\subfigure[IG6MP]{\includegraphics[width=0.48\textwidth]{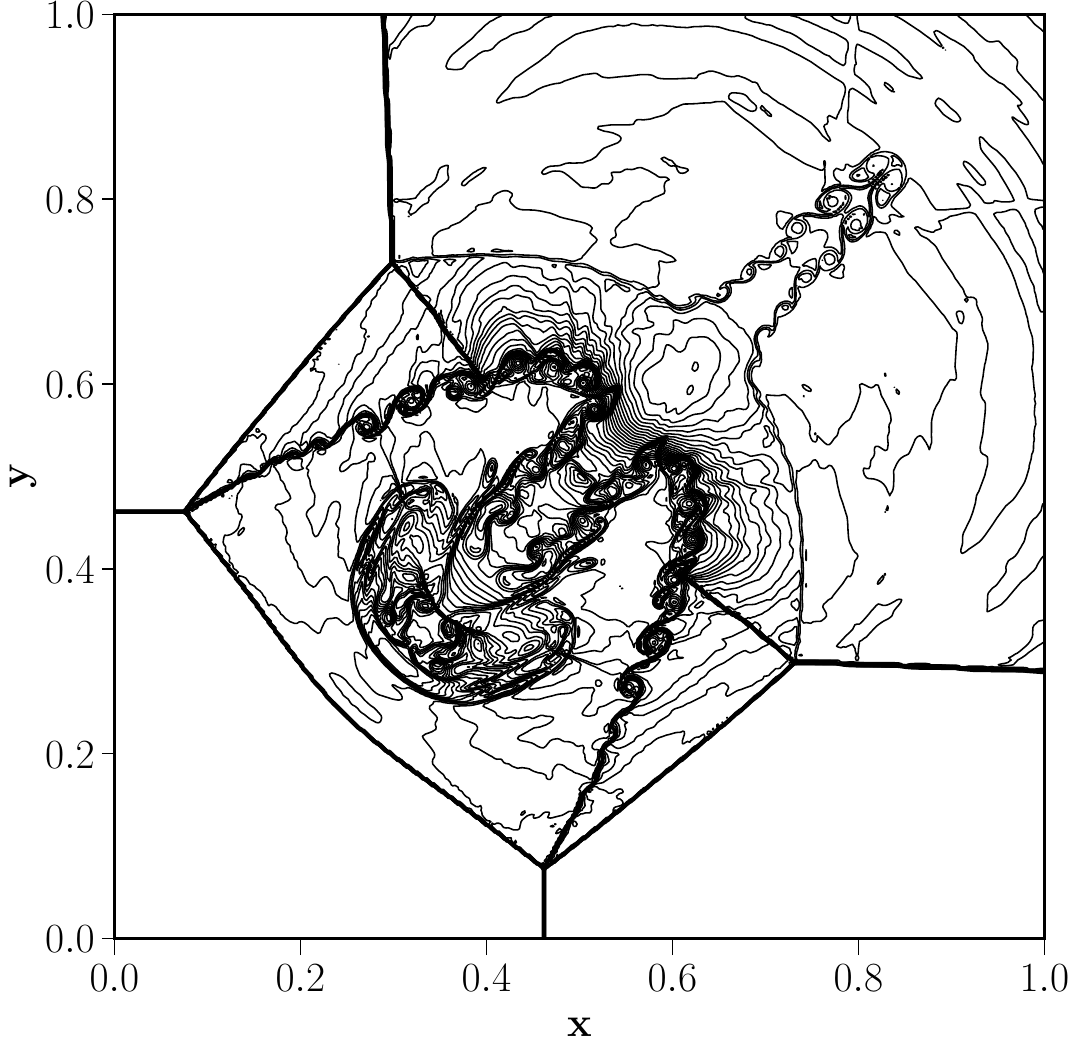}
\label{fig:IG6MC_RM}}
\subfigure[IG4MP]{\includegraphics[width=0.48\textwidth]{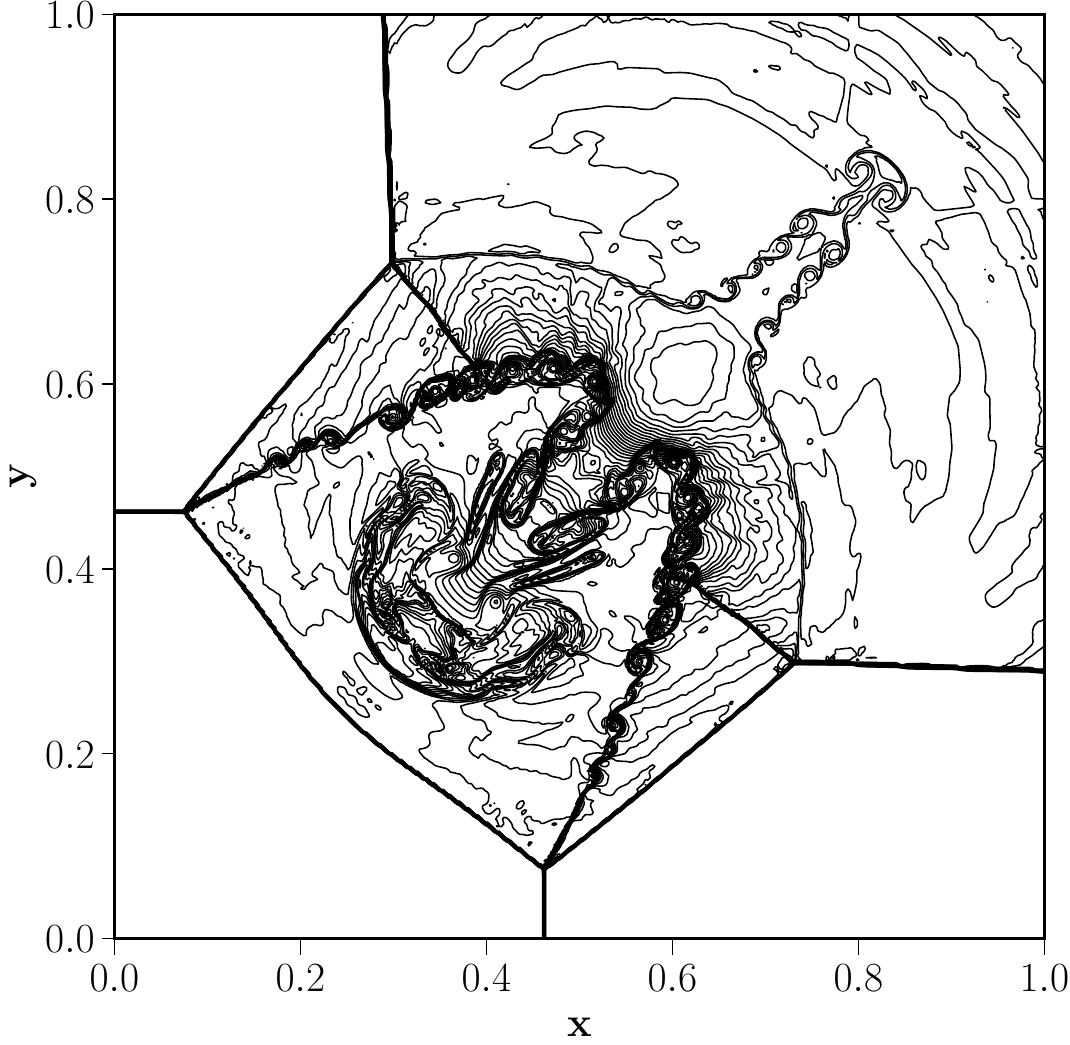}
\label{fig:IG4MP_RM}}
\caption{Density contours of the Riemann problem in Example \ref{ex:rp} for different schemes on a grid size of 400 $\times$ 400.}
\label{fig_riemann}
\end{onehalfspacing}
\end{figure}

The computed density contours are presented in Fig. \ref{fig_riemann}. The proposed schemes resolve better rollup behaviour along the slip lines compared to the \textcolor{black}{TENO5} scheme. The vortices indicate low numerical dissipation features for the proposed schemes.

\begin{example}\label{ex:rt} {{\color{black}Rayleigh-Taylor instability}}
\end{example}

The Rayleigh-Taylor instability occurs at the interface between fluids with different densities when acceleration is directed from the denser fluid to the less dense fluid. In this test case, two initial gas layers with different densities are subjected to unit magnitude's gravity, where the resulting acceleration is directed towards the less dense fluid. A small disturbance of the contact line triggers the instability. This problem has been extensively studied using high order shock-capturing schemes in the literature, see, e.g. \cite{Shi2002}, with the following initial conditions,

 \begin{equation}
\begin{aligned}
(\rho,u,v,p)=
\begin{cases}
(2.0,\ 0,\ -0.025\sqrt{\frac{5p}{3\rho}\cos(8\pi x)},\ 2y+1.00),&\quad 0\leq y< 0.5,\\
(1.0,\ 0,\ -0.025\sqrt{\frac{5p}{3\rho}\cos(8\pi x)},\ 1y+3/2),&\quad 0.5\leq y\leq 1.0,
\end{cases}
\end{aligned}
\label{eu2D_RT}
\end{equation}
over the computational domain $[0, 1/4]\times [0,1]$. Reflective boundary conditions are imposed on the right and left boundaries via ghost cells. The flow conditions are set to $\rho=1$, $p=2.5$, and $u=v=$0 on the top boundary and $\rho=2$, $p=1.0$, and $u=v=0$ on the bottom boundary with the specific heat ratio, $\gamma = 5/3$. The source term $S=(0,0,\rho, \rho v)$ is added to the Euler equations. We performed simulations on a uniform mesh of resolution $120 \times 480$ and the computations are conducted until $t = 1.95$.  

\begin{figure}[H]
\begin{onehalfspacing}
\centering\offinterlineskip
\subfigure[TENO5]{%
\includegraphics[width=0.20\textwidth]{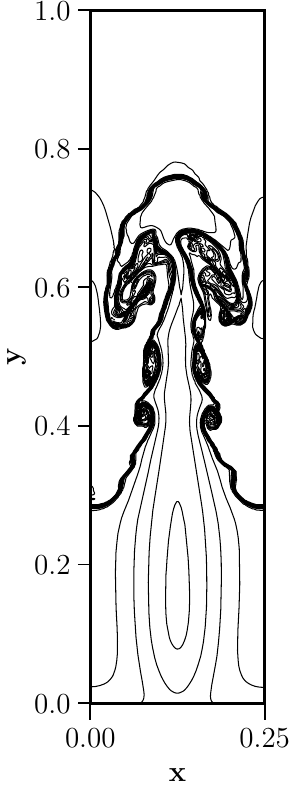}
\label{fig:RT_MP51}}
\subfigure[IG6MP]{%
\includegraphics[width=0.20\textwidth]{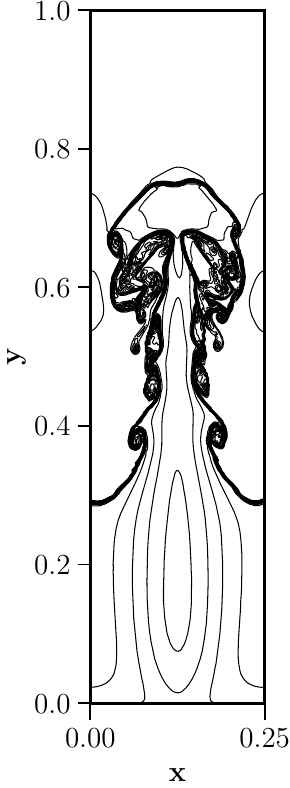}
\label{fig:RT_IG6T}}
\subfigure[IG4MP]{%
\includegraphics[width=0.20\textwidth]{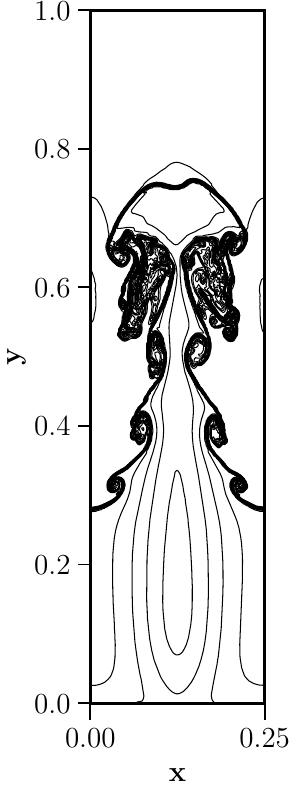}
\label{fig:RT_IG4T}}

\caption{Comparison of density contours obtained by different numerical schemes for the test case in Example \ref{ex:rt} on a grid size of 480 $\times$ 120.}
\label{fig:2d-RT1}
\end{onehalfspacing}
\end{figure}
Fig. \ref{fig:2d-RT1} shows the density distribution of the Rayleigh-Taylor instability problem. We can observe that the implicit gradient schemes produced more small vortices in the shear layer, indicating that they have better resolution to capture small scale features of the flow.

\begin{example}\label{ex:dmr}{Double Mach Reflection}
\end{example}

Next, we consider the double-Mach reflection problem proposed by \cite{woodward1984numerical}. In this test case, an unsteady planar shock-wave of Mach 10 impinges on an inclined surface of 30 degrees to the horizontal axis. This inclined surface is simplified by tilting the shock-wave to avoid modelling the oblique physical wall boundary. The near-wall jet structure and the vortex structures appearing from the contact discontinuity that emerges from the triple-point indicate the proposed scheme's numerical dissipation. Post-shock flow conditions are set at the left boundary, and zero gradient conditions are applied at the right boundary. At the bottom boundary, reflecting boundary conditions are applied in $x \in [1/6,3]$ and the post-shock conditions in $x \in \left[ 0, 1/6 \right]$. Furthermore, the exact solution of the moving shock is imposed at the upper at $y=1$ and is time-dependent. The computational domain is taken as $x \in [0,3], y \in [0,1]$ and the simulation is performed until $t=0.2$ on a grid of $768 \times 256$ cells.

\begin{equation}
\begin{aligned}
(\rho,u,v,p)=
\begin{cases}
&(8,\ 8.25 \cos 30^\circ,\ -8.25 \sin 30^\circ,\
116.5),\quad x<1/6+\frac{y}{\tan 60^\circ},\\
&(1.4,\ 0,\ 0,\
1),\quad\quad\quad\quad\quad\quad\quad\quad\quad\quad\quad\ x>
1/6+\frac{y}{\tan 60^\circ}.
\end{cases}
\end{aligned}
\label{eu2D_mach}
\end{equation}

\begin{figure}[H]
\centering\offinterlineskip
\subfigure[TENO5]{\includegraphics[width=0.3\textheight]{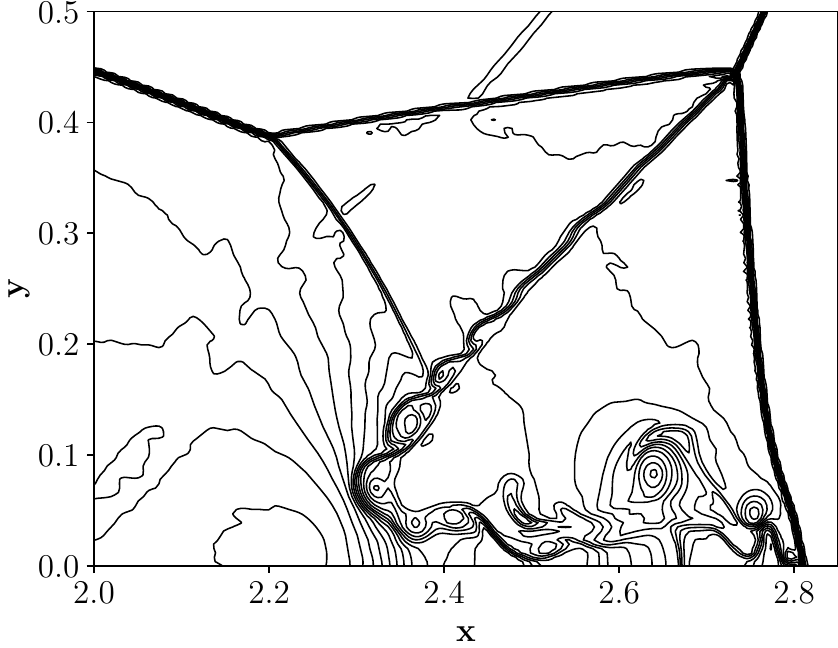}
\label{fig:MP5_DB}}
\subfigure[IG6MP]{\includegraphics[width=0.3\textheight]{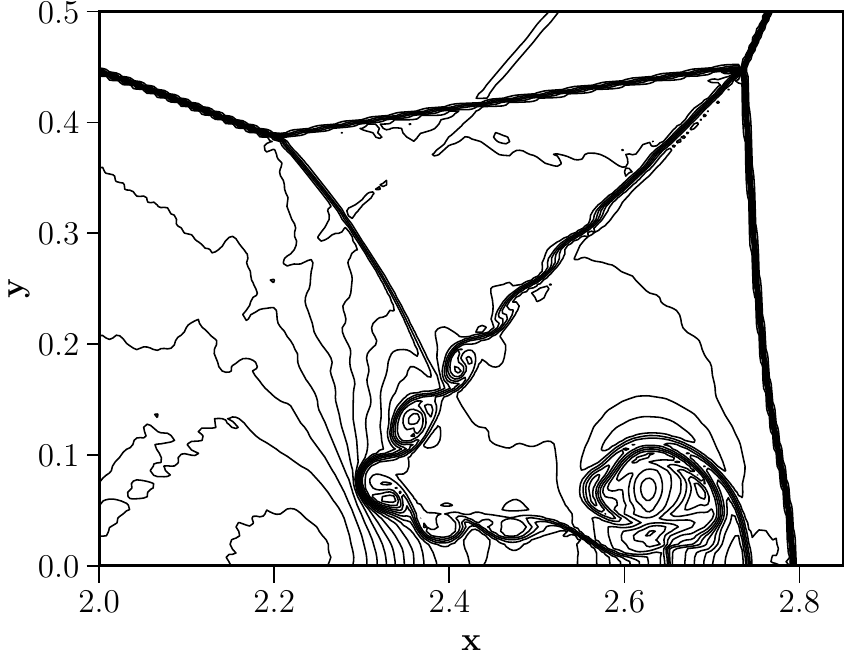}
\label{fig:HOCUS6_DB}}
\subfigure[IG4MP]{\includegraphics[width=0.3\textheight]{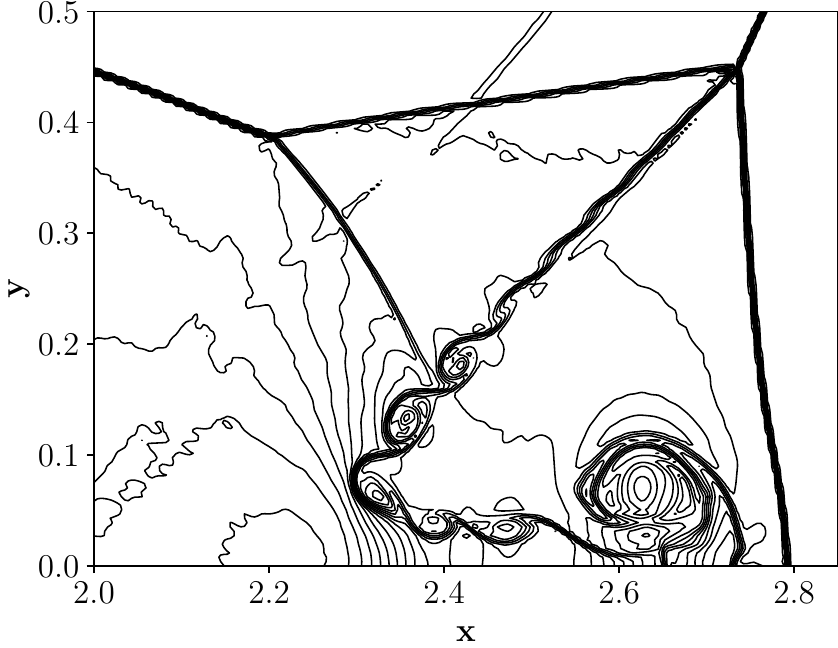}
\label{fig:MP5_DB1}}
\caption{Density contours in the blown-up region around the Mach stem for Example \ref{ex:dmr} on a grid size of 768 $\times$ 256.}
\label{fig_doublemach}
\end{figure}

Observations made from Figs. \ref{fig_doublemach} indicate that the IG4MP and the IG6MP schemes have better resolution of the Kelvin-Helmholtz (KH) instabilities \textcolor{black}{ than the TENO5 scheme}. Notably the resolution of the shear layers along the slip lines and the near-wall jet region are well resolved. It can be noted that the current IG4MP scheme is slightly better in resolving the shear layer along the slip line than the IG6MP scheme. 


\begin{example}\label{ex:TGV}{Inviscid Taylor-Green Vortex}
\end{example}

Next, we consider the three-dimensional inviscid Taylor-Green vortex problem, with initial conditions given by equation \eqref{itgv}. We use a domain of size $x,y,z \in [0,2\pi)$. Periodic boundary conditions are applied for all boundaries. The ratio of the specific heats of the gas is taken as $\gamma=5/3$. The simulations are performed until $t=10$ on a grid size of $64 \times 64 \times 64$. 

\begin{equation}\label{itgv}
\begin{pmatrix}
\rho \\
u \\
v \\
w \\
p \\
\end{pmatrix}
=
\begin{pmatrix}
1 \\
\sin{x} \cos{y} \cos{z} \\
-\cos{x} \sin{y} \cos{z} \\
0 \\
100 + \frac{\left( \cos{(2z)} + 2 \right) \left( \cos{(2x)} + \cos{(2y)} \right) - 2}{16}
\end{pmatrix}.
\end{equation}

This flow problem is essentially incompressible as the mean pressure is chosen to be very large. The Taylor-Green vortex is the simplest problem for analyzing the nonlinear transfer of kinetic energy among the different scales of the flow. It contains several physical processes that are key to understanding turbulence. The vortices in the initial flow stretch and produce smaller-scale features with time. This problem can be used as a test to examine the scale-separation ability of different schemes to under-resolved flow. We compare the ability of different schemes to preserve kinetic energy and also the growth of enstrophy in time, i.e., the sum of vorticity of all the vortex structures, indicating the schemes ability to preserve as many structures as possible. The enstrophy can be described as the integral of the square of the vorticity that can be computed as the integral of the magnitude of vorticity, $\overrightarrow{\omega}$, over the whole domain, 

\begin{equation}
Enstrophy = \sum_{cells}{} \left \| \overrightarrow{\omega} \right \|
\end{equation}

\begin{figure}[H]
\centering
\subfigure[Kinetic energy]{\includegraphics[width=0.48\textwidth]{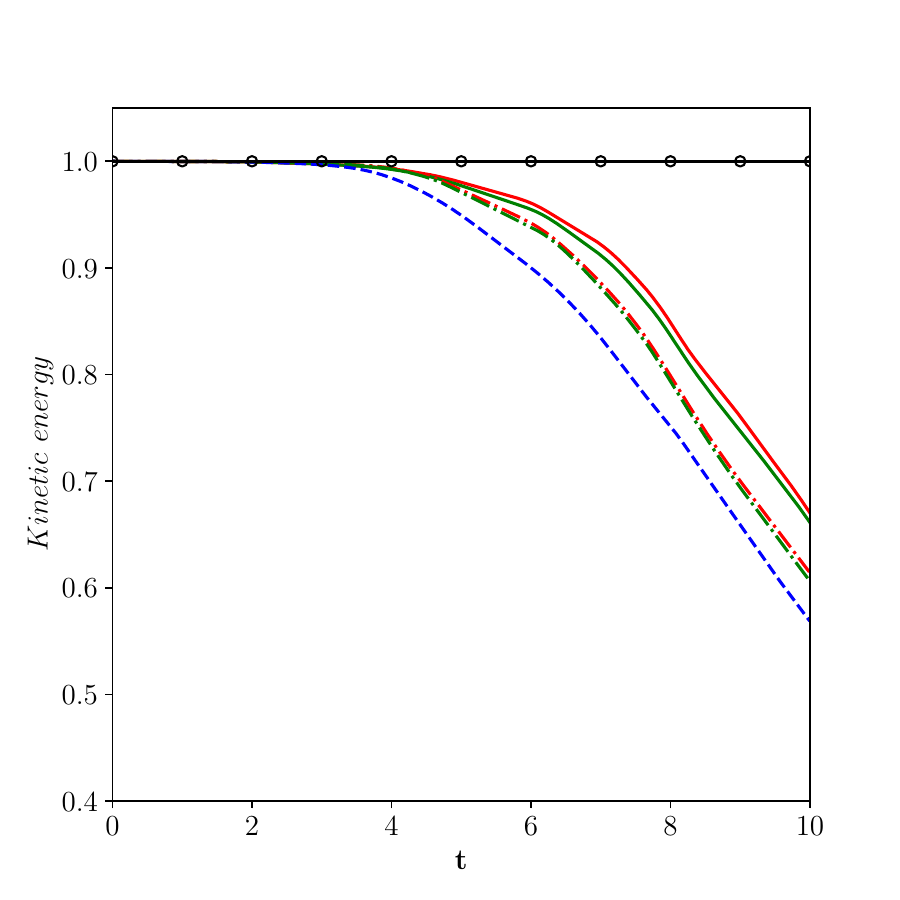}
\label{fig:TGV_KE}}
\subfigure[Enstrophy]{\includegraphics[width=0.48\textwidth]{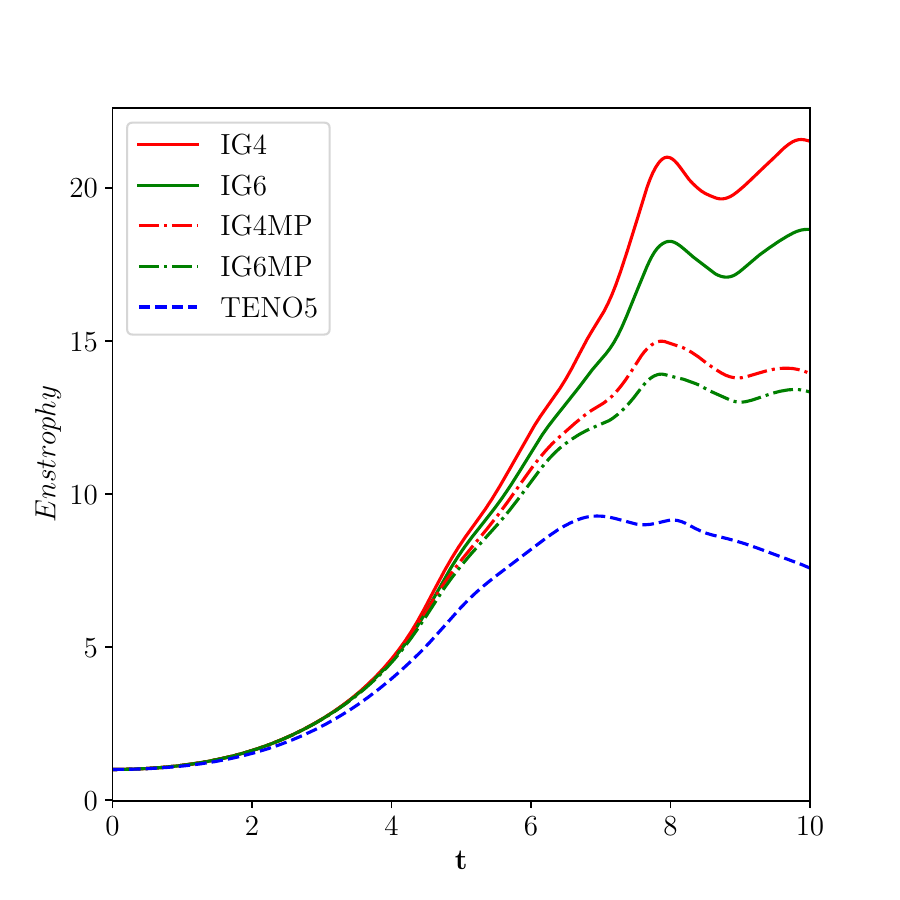}
\label{fig:TGV_ens}}
\caption{Normalised kinetic energy and enstrophy for different schemes presented in Example \ref{ex:TGV} on grid size of $64^3$. Solid line with circles: exact solution; solid red line: IG4; solid green line: IG6; dashed red line: IG4MP; dashed green line: IG6MP; dashed blue line: TENO5.}
\label{fig_TGV}
\end{figure}

Fig. \ref{fig_TGV} shows the normalised kinetic energy and normalised enstrophy with respect to the initial values for different schemes. \textcolor{black}{We have the following observations:}

\begin{itemize}
\item \textcolor{black}{First, we consider the kinetic energy evolution of all the numerical schemes considered here, shown in Fig. \ref{fig:TGV_KE}. The unlimited linear schemes IG4 and IG6 better preserved the kinetic energy than the limited nonlinear schemes. The nonlinear IGMP schemes preserved the kinetic energy better than the TENO5 scheme, with IG4MP slightly better than IG6MP.}
\item Next, we consider the enstrophy plot shown in Fig. \ref{fig:TGV_ens}. It can be observed that the present schemes outperform the TENO5 scheme significantly. Even though the BVD algorithm effectively captures discontinuities, there is still a significant difference between the kinetic energy and enstrophy values computed by the linear IG and nonlinear IGMP schemes, which can be improved in the future.  The statistics obtained by using different derivative schemes are presented in Appendix F of \textcolor{black}{Subramanium et al.} \cite{subramaniam2019high} by Subramaniam et al.. They have obtained the highest enstrophy by using spectral derivatives, and by using lower-order derivatives, they captured much lower enstrophy. They also noted that compact derivatives gave better results than explicit derivatives. Our results also indicate the same. An important advantage of the present IG schemes is the re-use of the velocity gradients in computing the enstrophy, which is the same used for the reconstruction of the interface states given by Equations (\ref{eq:cd4}), for IG4 and IG6 schemes, respectively. TENO5 does not have that advantage as it has to be computed separately. \textcolor{black}{For the IG4 and IG6 schemes, the native ``implicit gradient'' scheme that is used to compute the velocity gradients (as in Equation \ref{eqn:IG}) is used to compute the enstrophy as well, and it may lead to improved enstrophy values for the IG4 scheme.}

\end{itemize}

\subsection{Multi-dimensional test cases for Compressible Navier-Stokes equations}

\begin{example}\label{ex:vs}{Viscous Shock tube}
\end{example}

\begin{figure}[H]
\centering\offinterlineskip
\subfigure[IG6MP, $t$ = 1]{\includegraphics[width=0.3\textheight]{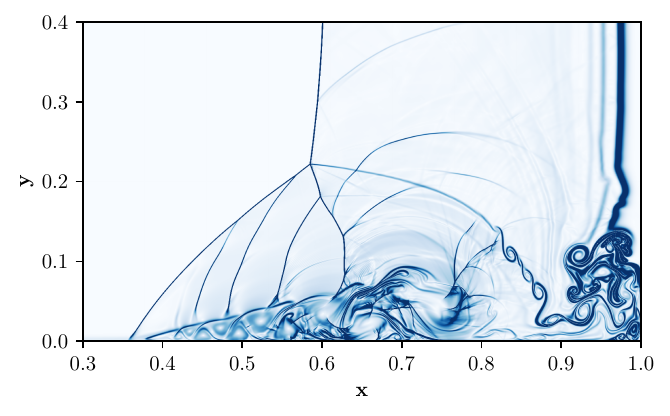}
\label{fig:IG6MP_VST_2500tp}}
\subfigure[IG4MP, $t$ = 1]{\includegraphics[width=0.3\textheight]{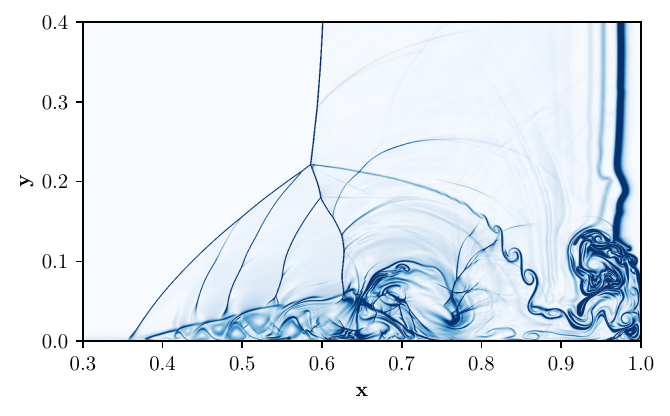}
\label{fig:IG4MP_VST_2500tp}}
\subfigure[Wall density profiles]{\includegraphics[width=0.25\textheight]{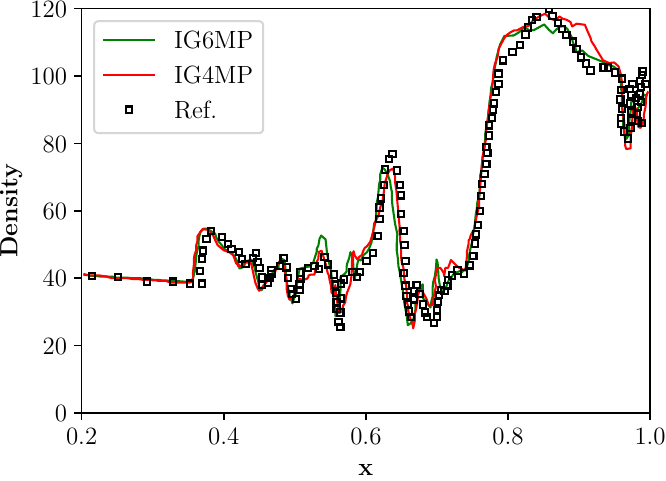}
\label{fig:wal}}
\caption{Density gradient plots using IG6MP and IG4MP, Example \ref{ex:vs}, for $Re=2500$ on a grid size of 2000 $\times$ 1000 and the wall density profile.}
\label{fig_VST-2500}
\end{figure}

We demonstrate that the proposed schemes can produce superior results even in viscous problems.  The viscous shock-tube problem of Daru and Tenaud \cite{daru2009numerical} is considered here. In this problem, the propagation of the shock wave and contact discontinuity leads to the developing of a thin boundary layer at the bottom wall. The shock wave interacts with this boundary layer once it reflects at the right wall. These interactions result in a complex vortex system, separation region, and a typical lambda-shaped shock pattern making it an ideal test case for evaluating high-resolution schemes. The initial conditions are:
\begin{equation}\label{vst}
\begin{aligned}
\left( {\rho , u,v, p} \right) = \left\{ \begin{array}{l}
\left( {120, 0 ,0,120/\gamma } \right),  \quad 0 < x < 0.5,\\
\left( {1.2 , 0 ,0, 1.2/\gamma } \right),  \quad 0.5 \le x < 1.
\end{array} \right.
\end{aligned}
\end{equation}
The domain for this test case is taken as $x \in [0,1], y \in [0,0.5]$. The initial conditions are given by Equation \eqref{vst}, with the ratio of specific heats of $\gamma = 7/5$. The flow is simulated for time $t=1$, keeping the Mach number of the shock wave at 2.37.  The problem is solved using the proposed schemes for a Reynolds number of $Re = 2500$ on a grid size of 2000 $\times$ 1000. The flow structures are much more complicated for this Reynolds number as the boundary layer separates at several points, giving rise to the development of highly complex vortex structures and the interactions between vortices and shock waves. Kundu et al. \cite{kundu2021investigation} carried out a fine grid simulation on 109 million cells for this Reynolds number evaluating the shear layer instabilities and vortex generation during the shock-wave boundary layer interaction. The numerical results obtained by the IGMP schemes are very close to their converged results (see Fig. 5 in \cite{kundu2021investigation}) despite using a grid 55 times smaller. Density distributions along the bottom wall, shown in Fig. \ref{fig:wal} also agree very well with those of Kundu et al. \cite{kundu2021investigation}. These results indicate that the proposed schemes can compute the multi-scale flows with shock waves in high resolution, even on coarse grids. The TENO5 scheme failed for this test case, and therefore, the results are not presented. Additionally, one advantage of the present implicit gradient method is sharing the velocity gradients between inviscid and viscous fluxes, which would not be possible with the TENO5 scheme.

\section{\textcolor{black}{Possible limitations of BVD approach}}\label{sec-5}

\subsection{\textcolor{black}{Unlimited schemes and analysis of BVD algorithm for ``low dissipation'' scheme.}}

\textcolor{black}{In this section, we present the difficulties of the BVD algorithm when a low-dissipation linear scheme is used. In this regard, the low dissipation IG4H (H stands for Hermite) considered by Chamarthi in \cite{CHAMARTHI2022105706} is used as the linear scheme as opposed to the IG4 approach considered earlier. In the IG4H approach, the second derivatives are computed as a function of the primitive variables and the gradients at the cell centers and are given by the following formula.}

\textcolor{black}{\begin{equation}\label{eqn:adam}
\mathbf{ {U}''}_{j}=\left(\frac{2}{ \Delta x^{2}}\right)\left(\mathbf{\hat {U}}_{j+1}-2\mathbf{\hat {U}}_{j}+\mathbf{\hat {U}}_{j-1}\right)-(\frac{1}{ 2 \Delta x})\left(\mathbf{ {U}'}_{j+1}-\mathbf{ {U}'}_{j-1}\right)
\end{equation}}

\textcolor{black}{Such an evaluation of the second derivatives will lead to a low dissipation scheme. In the equations (\ref{eqn:left_right}), we compute the first derivatives by using the optimized fourth-order compact derivatives given by Equation (\ref{eq:cd4}). Then, we use the explicit scheme given by Equation (\ref{eqn:adam}) to compute second derivatives. The spectral properties of the IG4H are already analyzed in \cite{CHAMARTHI2022105706}, and the readers can refer to the original paper.}

\begin{figure}[H]
\centering
\subfigure[Dispersion]{\includegraphics[width=0.45\textwidth]{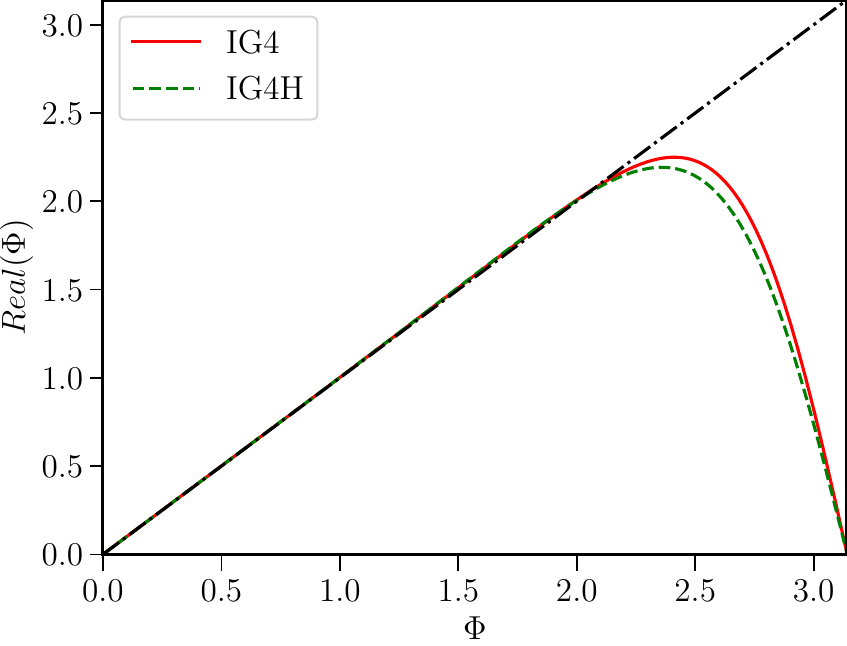}
\label{fig:dispersion_explicit}}
\subfigure[Dissipation]{\includegraphics[width=0.45\textwidth]{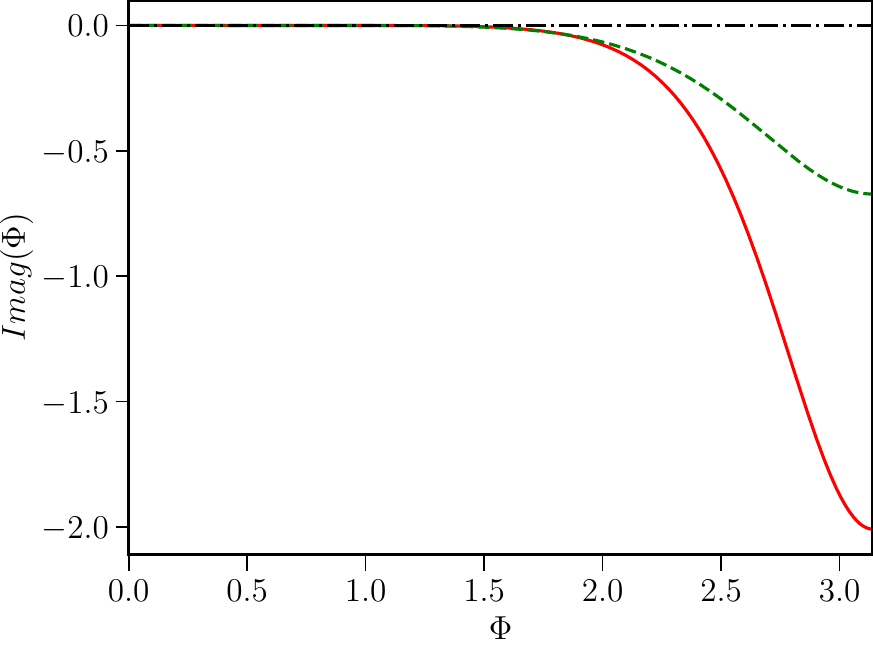}
\label{fig:dissipation_explicit}}
\caption{\textcolor{black}{Dispersion and Dissipation properties of the linear upwind schemes IG4 and IG4H.}}
\label{fig_disp_explicit}
\end{figure}

\textcolor{black}{One can see from Fig. \ref{fig_disp_explicit} that the IG4H scheme has lower dissipation than the IG4. Once again, we carried out the numerical simulations for the inviscid Taylor-Green vortex for the low dissipation scheme, Example \ref{ex:TGV}. Similar to the observations in Fig. \ref{fig_TGV}, the unlimited linear scheme IG4H gave better results than the IG4 and IG6 for kinetic energy and enstrophy, which also agrees with our theoretical observations from the Fourier analysis in Fig. \ref{fig_disp_explicit}.} 

\begin{figure}[H]
\centering
\subfigure[Kinetic energy]{\includegraphics[width=0.41\textwidth]{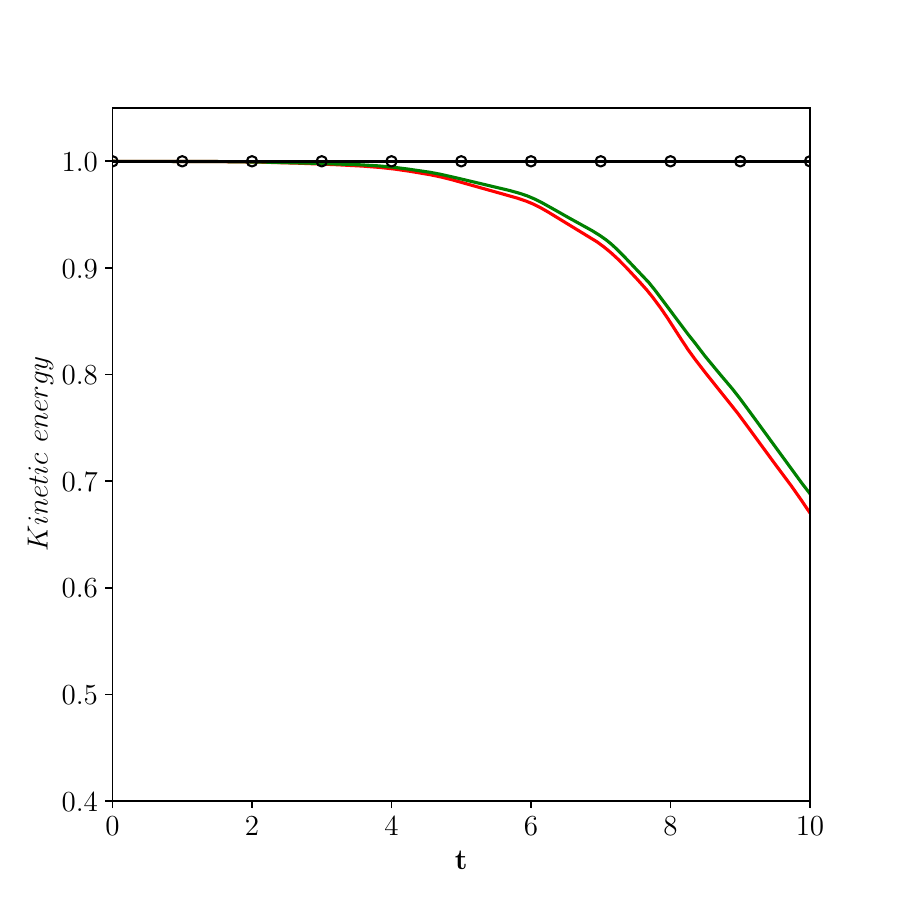}
\label{fig:TGV_KE_explicit}}
\subfigure[Enstrophy]{\includegraphics[width=0.4\textwidth]{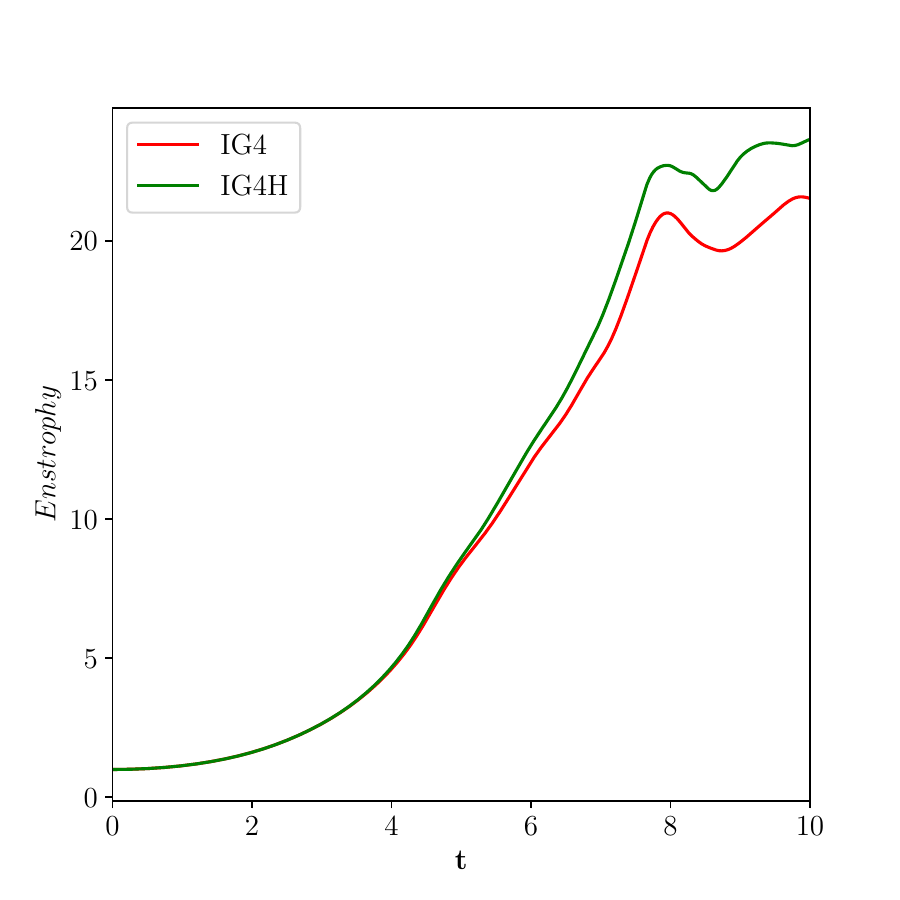}
\label{fig:TGV_ens_explicit}}
\caption{\textcolor{black}{Normalised kinetic energy and enstrophy for different schemes presented in Example \ref{ex:TGV} on grid size of $64^3$. Solid line with circles: exact solution; red line: IG4; green line: IG4H.}}
\label{fig_TGV_explicit}
\end{figure}

\textcolor{black}{The Taylor-Green vortex test case does not have discontinuities in the flow; therefore, the linear IG4H scheme may not encounter any difficulties. Next, we will consider test cases with discontinuities. Chamarthi and Frankel already performed such an analysis with a linear scheme in \cite{chamarthi2021high}. It is once again carried out with the IG4H and IG4 schemes here. Readers can see Fig. 14 and the corresponding analysis in  \cite{chamarthi2021high}.}\\

\textcolor{black}{Numerical experiments are carried out for the Examples \ref{sod} and \ref{Shu-Osher} using the linear schemes IG4H and IG4, i.e., without any nonlinear shock-capturing mechanism and compared with the exact solutions for better understanding of the BVD algorithm. All the numerical simulations are performed with a CFL of $0.2$. \textcolor{black}{Primitive} variables are interpolated to the cell interfaces for both IG4H and IG4 schemes. For the Sod and Lax problems in Example \ref{sod}, we considered $N=200$, and for the Shu-Osher problem in Example \ref{Shu-Osher}, we carried out the simulations on the grid size of $N=300$ and the following observations are made:}

\begin{figure}[H]
\centering\offinterlineskip
\subfigure[]{\includegraphics[width=0.38\textwidth]{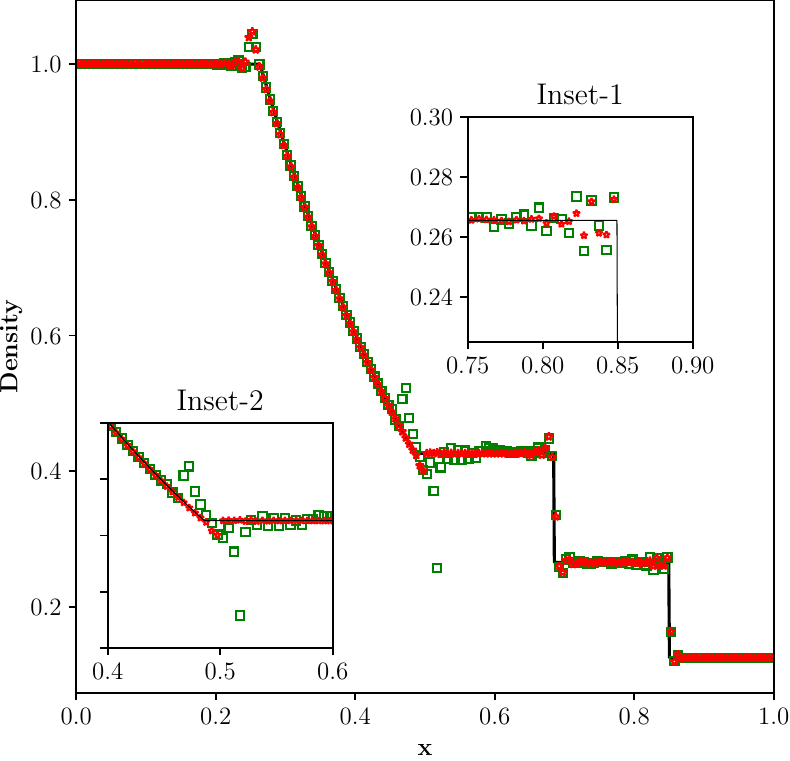}
\label{fig:unlim-sod}}
\subfigure[]{\includegraphics[width=0.38\textwidth]{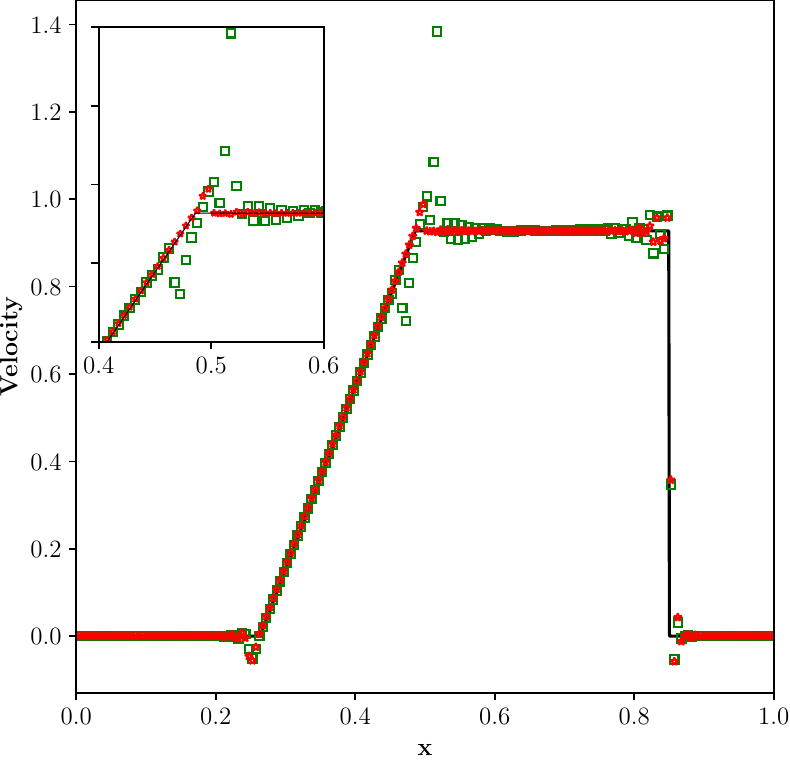}
\label{fig:unlim-sod_vel}}
\subfigure[]{\includegraphics[width=0.38\textwidth]{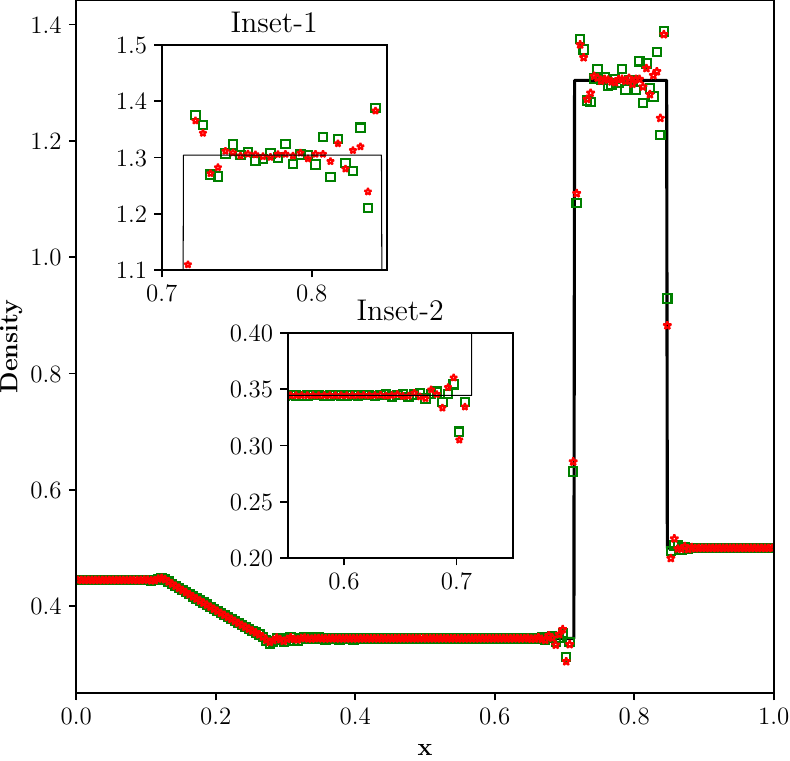}
\label{fig:unlim-lax_den}}
\subfigure[]{\includegraphics[width=0.38\textwidth]{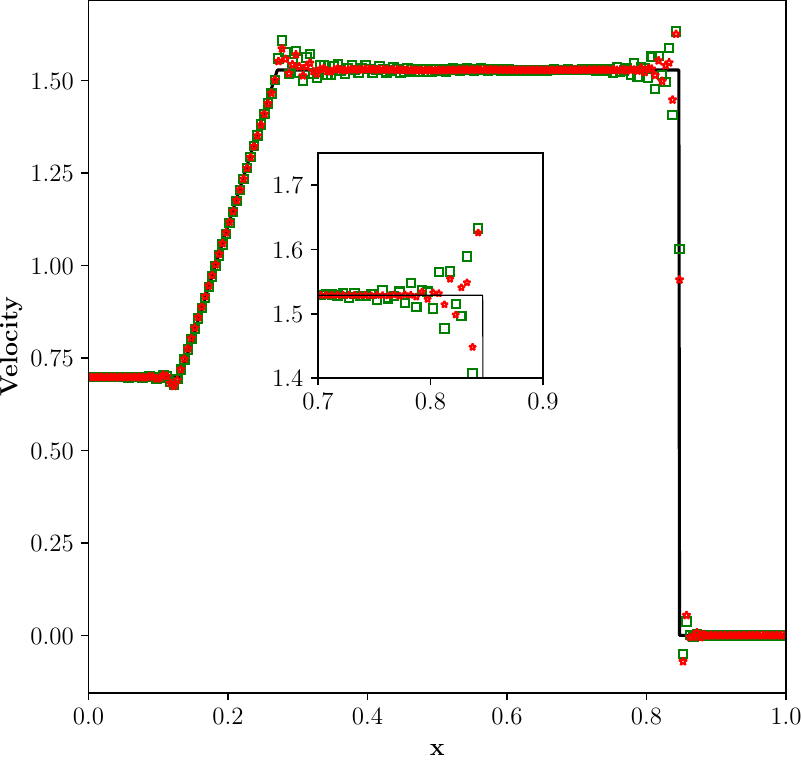}
\label{fig:unlim-lax}}
\subfigure[]{\includegraphics[width=0.38\textwidth]{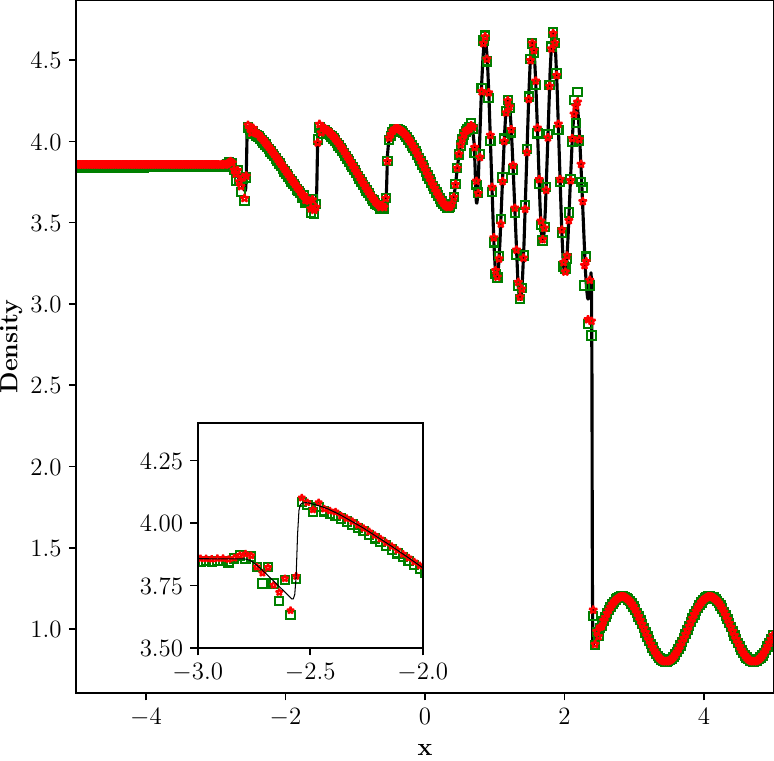}
\label{fig:unlim-shu}}
\subfigure[]{\includegraphics[width=0.38\textwidth]{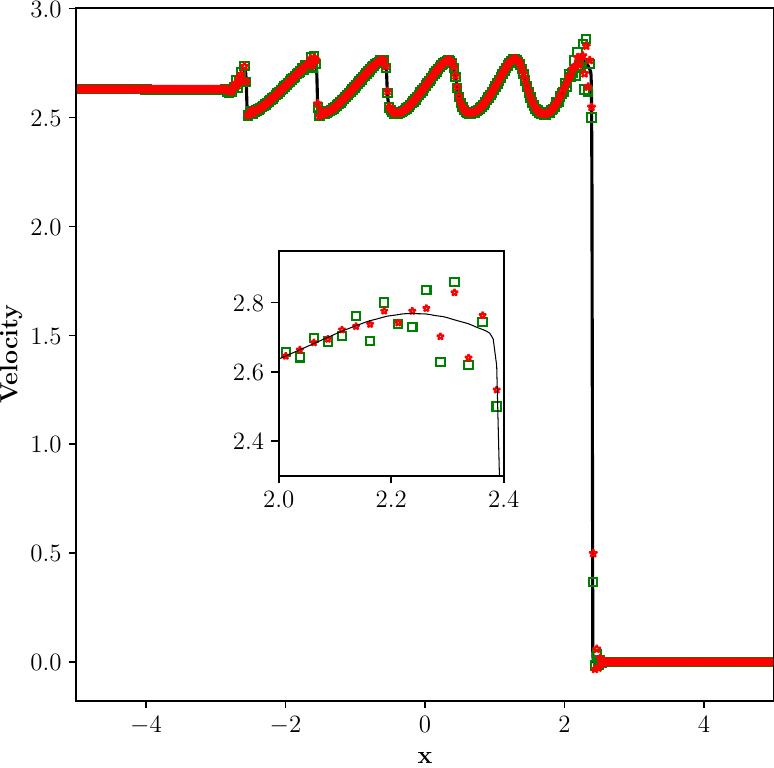}
\label{fig:unlim-shu-vel}}
\caption{Density and velocity profiles by IG4 and IG4H schemes for various test cases. Figs. (a) and (b) corresponds to the Sod problem, Figs. (c) and (d) corresponds to the Lax problem in Example \ref{sod},and Fig. (e) and (f) corresponds to the Shu-Osher problem, Example \ref{Shu-Osher}. Red stars: IG4 and Green squares: IG4H.}
\label{4E_unlim}
\end{figure}

\begin{itemize}
\item \textcolor{black}{Figs. \ref{fig:unlim-sod} and \ref{fig:unlim-sod_vel} show the density and velocity profiles for the Sod problem. We can see that the oscillations are only observed near the discontinuities for the IG4 scheme, whereas for the IG4H scheme, the oscillations are much more spread out and with higher amplitudes. Near the expansion wave, shown in inset-2 of Fig. \ref{fig:unlim-sod}, the IG4H scheme has more oscillations which are clearly due to the low dissipation property of the scheme in high wavenumber regions. Similarly, even in the velocity profile, the oscillations observed in IG4H are far more pronounced than in the IG4 scheme. Similar observations can also be made for the Lax problem from Figs. \ref{fig:unlim-lax_den} and \ref{fig:unlim-lax}. Therefore, using the BVD approach, it is sufficient to \textit{correct} fewer cells for the IG4 scheme, along with the Riemann solver, compared with the IG4H scheme.}
\item \textcolor{black}{Figs. \ref{fig:unlim-shu} and \ref{fig:unlim-shu-vel}  show the density and velocity profiles for the Shu-Osher problem using IG4 and IG4H schemes. Both the schemes have preserved the wave-like structures, but once again, the IG4H scheme has oscillations much more spread out in the velocity profile as shown in the inset of Fig.\ref{fig:unlim-shu-vel}.}
\end{itemize}

\textcolor{black}{With these observations, we have conducted further numerical experiments using the IG4H and MP5 schemes as a combination for the BVD algorithm. Similar to the approach used for the IG4 and IG6 schemes earlier, we compute the TBV of the IG4H scheme using the Equation (\ref{Eq:TBVIG4H}). \textcolor{black}{All the interface values at $j-\frac{3}{2}$, $j-\frac{1}{2}$, $j+\frac{1}{2}$ and $j+\frac{3}{2}$, for both $L$ and $R$} are modified according to the same algorithm as before, Equation (\ref{Eq:BVD-IG4H}), for each primitive variable. We denote the non-linear scheme as IG4H-BVD.}

\begin{equation}\label{Eq:TBVIG4H}
{TBV}_{j}^{IG4H}=\big|{U}_{j-\frac{1}{2}}^{L,IG4H}-{U}_{j-\frac{1}{2}}^{R,IG4H}\big|+\big|{U}_{j+\frac{1}{2}}^{L,IG4H}-{U}_{j+\frac{1}{2}}^{R,IG4H} \big|,
\end{equation} 

\begin{equation}\label{Eq:BVD-IG4H}
\text{if} \; {TBV}^{MP5} < {TBV}^{IG4H} \;  \left\{\begin{matrix}
{U}^{K,IG4H}_{j-\frac{3}{2}} = {U}^{{K},MP5}_{j-\frac{3}{2}}, \\ 
\\
{U}^{K,IG4H}_{j-\frac{1}{2}} = {U}^{{K},MP5}_{j-\frac{1}{2}}, \\ 
\\
{U}^{K,IG4H}_{j+\frac{1}{2}} = {U}^{{K},MP5}_{j-\frac{1}{2}}, \\ 
\\
{U}^{K,IG4H}_{j+\frac{3}{2}} = {U}^{{K},MP5}_{j+\frac{3}{2}}.
\end{matrix}\right.
\end{equation}

\textcolor{black}{To elucidate the difficulties concerning low dissipation IG4H scheme along with the BVD algorithm, we consider the blast wave test case of Woodward and Colella \cite{woodward1984numerical} with the following initial conditions:}

\begin{equation}
        \left( \rho,u,p \right) =
        \begin{cases}
            (1,0,1000), & \text{if } 0.0 \leq x < 0.1, \\
            (1,0,0.01), & \text{if } 0.1 \leq x < 0.8, \\
            (1,0,100), & \text{if } 0.8 \leq x \leq 1.0.
        \end{cases}\label{blast}
\end{equation}

\noindent \textcolor{black}{The case was solved on a computational domain $x = [0,1]$ with $N = 400$ uniformly distributed grid points until a final time, $t = 0.038$.}

\textcolor{black}{ One can see oscillations in the density profile for the IG4H scheme using the BVD algorithm. As it is observed from Fig. \ref{4E_unlim} that the IG4H scheme has oscillations much more spread out than the IG4 scheme, which indicates the number of cells that are to be corrected by the MP5 scheme will be significantly more. These oscillations are directly related to the low-dissipation property of the IG4H scheme.}

\textcolor{black}{ Similar \textit{oscillatory} results are observed with the P4T2 scheme of Deng et al. \cite{deng2019fifth}, see their Fig. 13, which is a combination of the fifth-order linear upwind and THINC schemes. Oscillations in the blast wave test case are also observed in all the follow-up papers of the P4T2 scheme, see Fig. 6 in  \cite{deng2020constructing} and oscillations are observed even in velocity profiles, see Fig. 13 of Ref. \cite{tann2020solution}. The fifth-order linear upwind scheme also has low dissipation properties similar to the IG4H scheme, which may contribute to such oscillations.}

\begin{figure}[H]
\centering\offinterlineskip
\includegraphics[width=0.26\textheight]{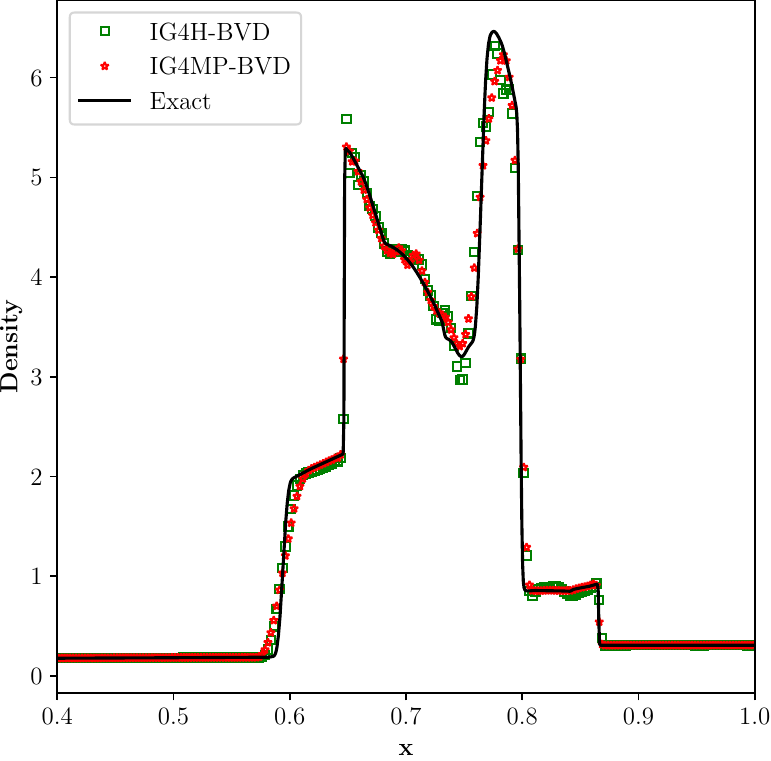}
\label{fig:blast_BVD}
\caption{\textcolor{black}{ Density profile for the blast wave test case, with initial conditions given by Equation (\ref{blast}), using the BVD approach for IG4H and IG4 schemes.}}
\label{fig_noBVD}
\end{figure}

\textcolor{black}{These results indicate that the BVD algorithm \textit{may encounter} difficulties with low dissipation schemes. To overcome the shortcomings of the BVD algorithm concerning the low dissipation IG4H scheme Chamarthi has proposed a different shock-capturing algorithm in \cite{CHAMARTHI2022105706}. The possible shortcoming of the BVD algorithm is that the no.of cells that needs to be corrected based on the TBV criteria in Equation (\ref{Eq:BVD-IG4H}) is challenging to judge. Chamarthi and Frankel \cite{chamarthi2021high} have studied various combinations of schemes with the BVD algorithm, and the readers can refer to the Appendix of \cite{chamarthi2021high}.}

\subsection{\textcolor{black}{High-order accuracy with \textit{linear} IG4 and IG6 schemes}}

\textcolor{black}{Fourier analysis is typically carried out for a linear advection equation \cite{deng2019fifth,van2021towards,chamarthi2022gradient}. Through Fourier analysis, it has been shown in the earlier subsection that the newly derived implicit gradient schemes are fourth-order accurate. However, when the proposed scheme is implemented for Euler equations, specifically for non-linear test cases, it was shown to be only second-order accurate in the Table \ref{tab:addlabel}. To understand this loss of accuracy, we carried out the following analysis. The introduction explains that the convective fluxes are computed using a Riemann solver, and the generic equation is rewritten below.}

\begin{equation}
\textcolor{black}{\mathbf{F}^{\rm Riemann}_{j+\frac{1}{2}}
= \frac{1}{2}
\left[
{\mathbf{F_L}}
+ 
{\mathbf{F_R}}
\right]
-
 \frac{1}{2} | {\mathbf{A}_{j+\frac{1}{2}}}|({\mathbf{Q}^R_{j+\frac{1}{2}}}-{\mathbf{Q}^L_{j+\frac{1}{2}}}),}
\label{eqn:Riemann_revisit}
\end{equation}

\textcolor{black}{In this paper, the interface fluxes $\mathbf{F_L}$ and $\mathbf{F_R}$ are computed from the $\mathbf{U_L}$ and $\mathbf{U_R}$ as the BVD algorithm only works with primitive variable reconstruction \cite{deng2019fifth,deng2020constructing,sun2016boundary,chamarthi2021high}. However, interface fluxes thus obtained will lead to loss of high-order accuracy and will only be second-order accurate \cite{van2021towards}, Equation (\ref{second-cry}).}

\begin{equation}\label{second-cry}
\textcolor{black}{\mathbf{F_L} \quad \text{computed as} \quad \mathbf{F(U_L)} \rightarrow \text{Second-order accurate}	}
\end{equation}

\textcolor{black}{If the interface fluxes are directly computed using the fluxes at the cell centres (denoted as \textbf{f}), one will attain the desired order of accuracy, Equation (\ref{high-happy}).}

\begin{equation}\label{high-happy}
\textcolor{black}{\mathbf{F_L} \quad \text{computed as} \quad  \mathbf{F(f_L)} \rightarrow \text{High-order accurate}}	
\end{equation}

\textcolor{black}{The \textit{linear} IG4 and IG6 schemes are reevaluated with the flux reconstruction as opposed to the primitive variable reconstruction in Equation \ref{eqn:IG}) and the results are shown in Table \ref{tab:high}. As a result, it can be seen that the genuine fourth-order accuracy for the Isentropic convecting vortex test case considered in Example \ref{euler-accuracy} has been obtained by both IG4 and IG6 schemes.}

{\begin{table}[H]
  \centering
  \caption{\textcolor{black}{High-order accuracy using flux reconstruction}}
    \begin{tabular}{crcrc}
          &       &       &       &  \\
    \hline
    N     & \multicolumn{1}{c}{IG4 (f)} &   Order    & \multicolumn{1}{c}{IG6 (f)} & Order \\
    \hline
    $25^2$ & \multicolumn{1}{c}{8.80E-04} & -     & \multicolumn{1}{c}{2.05E-03} & - \\
    \hline
    $50^2$ & \multicolumn{1}{c}{3.86E-05} & 4.51  & \multicolumn{1}{c}{6.32E-05} & 5.02 \\
    \hline
    $100^2$ & \multicolumn{1}{c}{2.02E-06} & 4.26  & \multicolumn{1}{c}{2.51E-06} & 4.65 \\
    \hline
    $200^2$ & 1.26E-07 & 4.01  & 1.35E-07 & 4.22 \\
    \hline
    \end{tabular}%
  \label{tab:high}%
\end{table}

\textcolor{black}{However, obtaining high-order accuracy concerning the nonlinear shock-capturing scheme is challenging when using the BVD algorithm. As explained earlier, the BVD algorithm works only with primitive variables and thus loses design order of accuracy. Therefore, obtaining both high-order accuracy and shock-capturing in the framework of the BVD algorithm is a challenging task (which is not presented in the open literature) and will be considered a future work.} \textcolor{black}{Alternatively, it may also be possible to use the approach used by Balsara and Kim \cite{balsara2016subluminal}, Zanotti and Dumbser \cite{zanotti2016efficient}, and Pidatella et al. \cite{pidatella2019semi} and still get high-order accuracy using primitive variable reconstruction with appropriate modifications in the current scheme for both inviscid and viscous fluxes.}


\section{Conclusions}\label{sec-6}

In this paper, {\color{black} linearly} high-order \textcolor{black}{implicit gradient schemes} have been developed {\color{black} based on a quadratic reconstruction combined with implicitly computed derivatives. Although the developed schemes are fourth-order accurate only for linear equations and the formal order of accuracy reduces to second-order for nonlinear equations, they still act as very low-dissipation/dispersion schemes and are capable of producing highly-resolved solutions even for problems with shock waves.} Important contributions and observations of the paper are summarized as follows

\begin{enumerate}
\item We proposed a novel approach of computing the cell interface values for \textcolor{black}{cell-centered conservative} framework where the gradients of reconstruction polynomials are computed by compact finite differences.  We have shown that fourth-order accuracy for linear problems can be achieved with a quadratic solution reconstruction if the derivatives are computed implicitly. 
\item We demonstrated that linearly high-order schemes can still serve as very low-dissipation/dispersion schemes for highly nonlinear problems with discontinuous solutions. Problem independent shock-capturing technique via the BVD algorithm gave superior results for several benchmark test cases involving shocks and small scale features.
\end{enumerate}

\section*{Acknowledgements}
A.S. and N. H are supported by Technion Fellowship during this work. A.S. dedicates this work to Vindhya Subrahmanyam, who passed away \textcolor{black}{in May 2021} (without her financial support, A.S. would not have traveled to Israel and carried \textcolor{black}{out any of his work}).

\section*{Appendix}
\renewcommand{\thesubsection}{\Alph{subsection}}
\textcolor{black} {\subsection{\textcolor{black}{Fifth-order monotonicity-preserving scheme}} \label{sec-appb}
In this Appendix the procedure for fifth-order monotonicity-preserving scheme of \cite{suresh1997accurate} is documented. For brevity, we only explain the procedure for the left interface values, $\mathbf{U}^{L,MP5}_{j+\half}$, since the right interface values, $\mathbf{U}^{R,MP5}_{j+\half}$ can be obtained via symmetry. The steps involved are as follows. }

\textcolor{black} {
\begin{equation}\label{eqn:mp5}
\mathbf{U}^{L,MP5}_{j+1/2}=\left\{\begin{array}{ll}
\mathbf{U}_{j+1 / 2}^{\text {Linear }} & \text { if }\left(\mathbf{U}_{j+1 / 2}^{\text {Linear }}-\mathbf{\hat{U}}_{j}\right)\left(\mathbf{U}_{j+1 / 2}^{\text {Linear }}-\mathbf{U}_{j+1 / 2}^{M P}\right) \leq 10^{-20}, \\
\mathbf{U}_{j+1 / 2}^{\text {Nonlinear }} & \text { otherwise, }
\end{array}\right.
\end{equation}
where
\begin{equation}\label{eqn:alpha}
\begin{aligned}
\mathbf{U}_{j+1 / 2}^{\text {Linear }} &=\frac{1}{60}(2 \mathbf{\hat{U}}_{j-2} - 13 \mathbf{\hat{U}}_{j-1} + 47 \mathbf{\hat{U}}_{j} + 27 \mathbf{\hat{U}}_{j+1} - 3 \mathbf{\hat{U}}_{j+2}),  \\
\mathbf{U}_{j+1 / 2}^{\text {Nonlinear }} &=\mathbf{U}_{j+1 / 2}^{\text {Linear }}+\operatorname{minmod}\left(\mathbf{U}_{j+1 / 2}^{\min }-\mathbf{U}_{j+1 / 2}^{\text {Linear }}, \mathbf{U}_{j+1 / 2}^{\max }-\mathbf{U}_{j+1 / 2}^{\text {Linear }}\right), \\
\mathbf{U}_{j+1 / 2}^{M P} &=\mathbf{U}_{j}+\operatorname{minmod}\left[\mathbf{\hat{U}}_{j+1}-\mathbf{\hat{U}}_{j}, 7\left(\mathbf{\hat{U}}_{j}-\mathbf{\hat{U}}_{j-1}\right)\right], \\
\mathbf{U}_{j+1 / 2}^{\min } &=\max \left[\min \left(\mathbf{\hat{U}}_{j}, \mathbf{\hat{U}}_{j+1}, \mathbf{U}_{j+1 / 2}^{M D}\right), \min \left(\mathbf{\hat{U}}_{j}, \mathbf{U}_{j+1 / 2}^{U L}, \mathbf{U}_{j+1 / 2}^{L C}\right)\right], \\
\mathbf{U}_{j+1 / 2}^{\max } &=\min \left[\max \left(\mathbf{\hat{U}}_{j}, \mathbf{\hat{U}}_{j+1}, \mathbf{U}_{j+1 / 2}^{M D}\right), \max \left(\mathbf{\hat{U}}_{j}, \mathbf{U}_{j+1 / 2}^{U L}, \mathbf{U}_{j+1 / 2}^{L C}\right)\right], \\
\mathbf{U}_{j+1 / 2}^{M D} &=\frac{1}{2}\left(\mathbf{\hat{U}}_{j}+\mathbf{\hat{U}}_{j+1}\right)-\frac{1}{2} d_{j+1 / 2}^{M}, \\
\mathbf{U}_{j+1 / 2}^{U L} &=\mathbf{\hat{U}}_{j}+4\left(\mathbf{\hat{U}}_{j}-\mathbf{\hat{U}}_{j-1}\right), \\
\mathbf{U}_{j+1 / 2}^{L C} &=\frac{1}{2}\left(3 \mathbf{\hat{U}}_{j}-\mathbf{\hat{U}}_{j-1}\right)+\frac{4}{3} d_{j-1 / 2}^{M}, \\
d_{j+1 / 2}^{M} &=\operatorname{minmod}\left(4 d_{j}-d_{j+1}, 4 d_{j+1}-d, d_{j},d_{j+1}\right), \\
d_{j} &=\mathbf{\hat{U}}_{j-1}-2 \mathbf{\hat{U}}_{j}+\mathbf{\hat{U}}_{j+1},
\end{aligned}
\end{equation}}

\textcolor{black} {where,
\begin{equation}
minmod(a,b) = \half \left ( sign(a)+sign(b) \right ) min(|a|,|b|), 
\end{equation}}

\textcolor{black} {\noindent \textcolor{black}{The MP5 reconstruction involves the transforming of primitive variables into characteristic variables and readers are referred to section 2.3} in \cite{chamarthi2021high}.}

\bibliographystyle{elsarticle-num}

\end{document}